\documentclass[12pt]{amsart}

\usepackage{amsmath,xspace,amssymb,mathrsfs}
\usepackage{color}

\input xy
\xyoption{all}
\xyoption{2cell}
\UseAllTwocells
\CompileMatrices

\newcommand{\colimit}{\operatorname{colim}}

\renewcommand{\phi}{\varphi}

\newcommand{\hm}{\operatorname{Hom}}
\newcommand{\Supp}{\operatorname{Supp}}

\newcommand{\Hom}{\operatorname{Hom}}
\newcommand{\Fib}{\operatorname{Fib}}

\newcommand{\ou}{\operatorname{Op}}
\newcommand{\Id}{\operatorname{Id}}
\newcommand{\sets}{\operatorname{\mathfrak{sets}}}

\newcommand{\Ab}{\operatorname{\mathfrak{Ab}}}
\newcommand{\can}{\operatorname{\mathfrak{c}}}
\newcommand{\psh}{\operatorname{\mathfrak{psh}}}
\newcommand{\sh}{\operatorname{\mathfrak{sh}}}
\newcommand{\leta}{\operatorname{a}}
\newcommand{\letb}{\operatorname{b}}
\newcommand{\letc}{\operatorname{c}}
\newcommand{\Top}{\operatorname{\mathfrak{Top}}}
\newcommand{\Ima}{\operatorname{Im}}

\newtheorem{proposition}{Proposition}[section]

\theoremstyle{definition}

\newtheorem{bk}[proposition]{}

\begin{document}

\title{Abstract sheaf theory}

\author{Abolfazl Tarizadeh}
\address{Department of Mathematics, Faculty of Basic Sciences, University of Maragheh \\
P. O. Box 55181-83111, Maragheh, Iran.
 }
\email{ebulfez1978@gmail.com}

\date{}

\begin{abstract} In this article, the theory of sheaves is studied from a categorical point of view. This perspective vastly generalizes the usual theory of sheaves of sets to a more abstract setting which allows us to investigate the theory of sheaves with values in an arbitrary category.
\end{abstract}

\maketitle

\section{Introduction}

This work is an introduction to the sheaf theory which the author wrote it down in the summer of 2014. We have tried to present the material in the way that an undergraduate student could be able to read them. Obviously, an advanced reader does not require to follow all of the material line-by-line while who can prove many proofs by himself/herself. In preparing this work the main sources which we have used are including \cite{Grothendieck} and \cite{Johan}. The titles of the sections should be sufficiently explanatory.\\

When we are studying algebraic geometry, we are practically involved with sheaves with values in various categories such as sets, groups, rings, topological rings, modules and etc. Thus studying the theory of sheaves from a more abstract (categorical) point of view gives us a vast insight into the theory and at the same time also gives us a comfortability to get ride of verifying a fact each time for various categories. Let us clarify the picture by giving an explicit example. Let $X$ be a topological space, $\mathfrak{B}$ a basis for the opens of $X$, $\mathscr{F}$ and $\mathscr{G}$ two sheaves over $X$ with values in a category $\emph{\textbf{K}}$. Also assume that there are morphisms $\phi_{U}:\mathscr{F}(U)\rightarrow\mathscr{G}(U)$ in $\emph{\textbf{K}}$ for all $U\in\mathfrak{B}$ which are compatible with the restrictions. Then there is a unique morphism of sheaves $\Phi:\mathscr{F}\rightarrow\mathscr{G}$ such that $\Phi_{U}=\phi_{U}$ for
all $U\in\mathfrak{B}$. This is proven easily when we studying the theory of sheaves from an abstract point of view and has the premier advantage that we do not need to verify this fact each time for the involved category in the context.

\section{Sheaves with values in a category}

\begin{bk}\label{I: 3.1.2} Let $X$ be a topological space, one can associate to $X$ a category $\ou(X)$, so-called the category of open subsets of $X$, whose objects are the open subsets of $X$ and morphisms are the inclusion maps. Every contravariant functor $\mathscr{F}$ from $\ou(X)$ into a category $\emph{\textbf{K}}$ is called a presheaf on $X$ with values in the category $\emph{\textbf{K}}$. Let $\mathscr{F}$ be a presheaf on $X$ with values in the category $\emph{\textbf{K}}$, therefore $\mathscr{F}$ attaches, to each open subspace $U$ of $X$ an object $\mathscr{F}(U)$ of $\emph{\textbf{K}}$ and
for each two open subspaces $U, V$ of $X$ with $U\subseteq V$, a morphism $\rho_{U}^{V}=\mathscr{F}(i_{U}^{V}):
\mathscr{F}(V)\rightarrow\mathscr{F}(U)$ in $\emph{\textbf{K}}$ ($i_{U}^{V}:U\rightarrow V$ is the inclusion map) which is usually called the restriction morphism of $\mathscr{F}(V)$ into $\mathscr{F}(U)$, so that $\rho_{U}^{U}$ is the identity morphism and if $U\subseteq V\subseteq W$ then $\rho_{U}^{W}=\rho_{U}^{V}\circ\rho_{V}^{W}$. \\

Suppose $\mathscr{F}$ and $\mathscr{G}$ are two presheaves on $X$ with values in a category $\emph{\textbf{K}}$, a morphism of presheaves $u:\mathscr{F}\rightarrow\mathscr{G}$ is just a morphism of functors from $\mathscr{F}$ into $\mathscr{G}$. Therefore, a morphism of presheaves $\mathscr{F}\rightarrow\mathscr{G}$ is an element $u=(u_{U})$ of $\prod\limits_{U\in\ou(X)}\Hom_{\emph{\textbf{K}}}(\mathscr{F}(U),  \mathscr{G}(U))$ such that for each two open subspaces $U,V$ of $X$ with $V\subseteq U$, the following diagram is commutative $$\xymatrix{
\mathscr{F}(U)\ar[r]^{u_{U}} \ar[d]^{\rho_{V}^{U}} & \mathscr{G}(U)\ar[d]^{\varrho_{V}^{U}} \\ \mathscr{F}(V)\ar[r]^{u_{V}} &\mathscr{G}(V).}$$
Hence the class of morphisms of presheaves from $\mathscr{F}$ into $\mathscr{G}$ is a set.\\

\begin{bk}\label{I: 3.1.1} A presheaf $\mathscr{F}$ over $X$ with values in a category $\emph{\textbf{K}}$ is called a $\emph{sheaf}$ if it satisfies in the following axiom:\\
($\textbf{F}$) For each open subspace $U$ of $X$, for each open covering $U=\bigcup\limits_{\alpha\in I}U_{\alpha}$ and for each object $T$ of  $\emph{\textbf{K}}$, the map $f\rightsquigarrow (\rho_{\alpha}\circ f)_{\alpha}$ is a bijection from $\hm_{\emph{\textbf{K}}}\big(T,\mathscr{F}(U)\big)$ onto the set of families $(f_{\alpha})\in\prod\limits_{\alpha\in I}\hm_{\emph{\textbf{K}}}\big(T,\mathscr{F}(U_{\alpha})\big)$ such that for each pair of indices $(\alpha,\beta)$, $\rho_{\alpha\beta}\circ f_{\alpha}=\rho_{\beta\alpha}\circ f_{\beta}$ where $\rho_{\alpha}:
\mathscr{F}(U)\rightarrow\mathscr{F}(U_{\alpha})$
and $\rho_{\alpha\beta}:
\mathscr{F}(U_{\alpha})\rightarrow\mathscr{F}(U_{\alpha}\cap U_{\beta})$ are the restriction morphisms. \\
\end{bk}

We call a category non-empty if its class of objects is non-empty. Denote by $\sh_{\emph{\textbf{K}}}(X)$ the category of sheaves over $X$ with values in $\emph{\textbf{K}}$, then $\sh_{\emph{\textbf{K}}}(X)$ is a non-empty category if and only if the category $\emph{\textbf{K}}$ has a terminal (final) object. Because first suppose that $\sh_{\emph{\textbf{K}}}(X)$ is a non-empty category and then take a sheaf $\mathscr{F}$ over $X$ with values in $\emph{\textbf{K}}$; consider the empty covering $(U_{\alpha})_{\alpha\in I}$ (i.e. $I=\emptyset$) for the empty set $U=\emptyset$, then for each object $T$ of $\emph{\textbf{K}}$, the set $\prod\limits_{\alpha\in I=\emptyset}\Hom_{\emph{\textbf{K}}}(T, \mathscr{F}(U_{\alpha}))=\{\emptyset\}$ is a singleton set whose its only element is the empty sequence and the map $$\Hom_{\emph{\textbf{K}}}(T, \mathscr{F}(\emptyset))\rightarrow\prod\limits_{\alpha\in I=\emptyset}\Hom_{\emph{\textbf{K}}}(T, \mathscr{F}(U_{\alpha}))=\{\emptyset\}$$ is bijective, therefore $\mathscr{F}(\emptyset)$ is a terminal object of $\emph{\textbf{K}}$. Conversely suppose that $D$ is a terminal object of $\emph{\textbf{K}}$, then $U\rightsquigarrow\mathscr{F}(U)=D$ is a sheaf over $X$ with values in $\emph{\textbf{K}}$. \\

A presheaf $\mathscr{F}$ on $X$ with values in the category $\emph{\textbf{K}}$ is a sheaf if and only if for each object $T$ of $\emph{\textbf{K}}$, the assignment $U\rightarrow\Hom_{\emph{\textbf{K}}}(T, \mathscr{F}(U))$ is a sheaf on $X$ with values in the category of sets. \\
\end{bk}

\begin{bk}\label{tanriverdi} Let $X$ be a topological space, $\emph{\textbf{K}}$ a category which admits products and $\mathscr{F}$ a presheaf over $X$ with values in $\emph{\textbf{K}}$. Then it is an interesting exercise to see that $\mathscr{F}$ is a sheaf if and only if for each open subset $U$ of $X$, for each open covering $U=\bigcup\limits_{\alpha\in I}U_{\alpha}$ the following diagram is an equalizer diagram
$$\mathscr{F}(U)\rightarrow\prod\limits_{\alpha\in I} \mathscr{F}(U_{\alpha})\rightrightarrows\prod\limits_{(\alpha,\beta)\in I^{2}} \mathscr{F}(U_{\alpha}\cap U_{\beta})$$ where the unnamed arrows are the unique morphisms which arise from the products.\\
\end{bk}

\begin{bk}\label{1.2} Let $\mathscr{F}$ be a presheaf on $X$ with values in the category of sets (in this case it is called a presheaf of sets on $X$). As a consequence of \ref{tanriverdi} we obtain that $\mathscr{F}$ is a sheaf if and only if for every open covering $(U_{\alpha})$ of an open subspace $U$ of $X$, it satisfies in the following conditions:\\
$(\mathbf{G1})$ If there exist some elements $s, t\in\mathscr{F}(U)$ so that $\rho_{\alpha}(s)=\rho_{\alpha}(t)$ for each $\alpha$, then $s=t$.\\
$(\mathbf{G2})$ If for each $\alpha$ there exists some $s_{\alpha}\in\mathscr{F}(U_{\alpha})$ so that for each pair of indices $(\alpha, \beta)$, $\rho_{\alpha\beta}(s_{\alpha})=\rho_{\beta\alpha}(s_{\beta})$ then there exists some element (in fact unique element according to $(\mathbf{G1})$) $s\in\mathscr{F}(U)$ so that $\rho_{\alpha}(s)=s_{\alpha}$ for all $\alpha$.\\
\end{bk}

\begin{bk}\label{I: 3.1.5} Let $\mathscr{F}$ be a presheaf (resp. a sheaf) on $X$ with values in a category $\emph{\textbf{K}}$, $U$ an open subspace of $X$, the functor
$V\rightarrow\mathscr{F}(V)$ for the open subspaces $V$ of $U$ constitute a presheaf (resp. a sheaf) on $U$ with values in $\emph{\textbf{K}}$, which is called the presheaf (resp. sheaf) induced via $\mathscr{F}$ on $U$ and it is denoted by $\mathscr{F}_{|_{U}}$.
For each morphism $u:\mathscr{F}\rightarrow\mathscr{G}$ of presheaves with values in $\emph{\textbf{K}}$, we denote by $u_{|_{U}}$ the morphism $\mathscr{F}_{|_{U}}\rightarrow\mathscr{G}_{|_{U}}$ induced via $u$ which is formed by $u_{V}$ for the open subspaces $V\subseteq U$.\\
\end{bk}

\begin{bk}\label{I: 3.1.6} Suppose that the category  $\emph{\textbf{K}}$ admits the inductive (direct) limits. Let $\mathscr{F}$ be a presheaf on $X$ with values in $\emph{\textbf{K}}$, the fiber of  $\mathscr{F}$ at the point $x\in X$, denoted by $\mathscr{F}_{x}$, is an object of $\emph{\textbf{K}}$ which is indeed the inductive limit of  $\mathscr{F}(U)$ and the morphisms $\rho_{V}^{U}: \mathscr{F}(U)\rightarrow\mathscr{F}(V)$ where $U$ and $V$ belong to the set of open neighborhoods of $x$ which is a filtered (directed) set ordered by inclusion (more precisely, for each two open neighborhoods $U,V$ of $x$ in $X$, we say that $U< V$ if $V\subset U$ ). Therefore $\mathscr{F}_{x}=\colimit_{x\in U}\mathscr{F}(U)$ where $U$ runs through the open neighborhoods of $x$. Moreover, if $u: \mathscr{F}\rightarrow\mathscr{G}$ be a morphism of presheaves with values in $\emph{\textbf{K}}$, for each point $x\in X$, we define $u_{x}:\mathscr{F}_{x}\rightarrow\mathscr{G}_{x}$ the unique morphism in $\emph{\textbf{K}}$ as the inductive limit of the morphisms $u_{U}: \mathscr{F}(U)\rightarrow\mathscr{G}(U)$ where $U$ runs through the open neighborhoods of $x$. For any such open subspace $U$ the following diagram is also commutative $$\xymatrix{
\mathscr{F}(U) \ar[r]^{u_{U}} \ar[d]^{} & \mathscr{G}(U) \ar[d]^{} \\ \mathscr{F}_{x}\ar[r]^{u_{x}} & \mathscr{G}_{x}} $$ where the vertical arrows are the canonical morphisms.
If $\mathscr{F}$ is a sheaf of sets on $X$, for each open subspace $U$ of $X$, the elements of $\mathscr{F}(U)$ are called the sections of $\mathscr{F}$ over $U$ and we usually write $\Gamma(U, \mathscr{F})$ instead of $\mathscr{F}(U)$. If $s\in\Gamma(U, \mathscr{F})$ and $V$ is an open subspace of $X$ contained in $U$ then we write $s|_{V}$ instead of $\rho_{V}^{U}(s)$. For each $x\in U$, the image of $s$ under the canonical morphism $\mathscr{F}(U)\rightarrow\mathscr{F}_{x}$ is called the \emph{germ} of $s$ at the point $x$ and it is denoted by $s_{x}$ if there is no confusion. Moreover if $u:\mathscr{F}\rightarrow\mathscr{G}$ is a morphism of sheaves of sets we shall often write $u(s)$ instead of $u_{U}(s)$ for all $s\in\Gamma(U, \mathscr{F})$. If $\mathscr{F}$ is a sheaf of abelian groups, or of rings, or of modules, then the set of points $x\in X$ so that $\mathscr{F}_{x}\neq\{0\}$ is called the support of  $\mathscr{F}$, denoted by $\Supp(\mathscr{F})$; this set is not necessarily closed in $X$.\\
\end{bk}

\section{Presheaves over a basis of opens}

Throughout this section, $\emph{\textbf{K}}$ is a category which admits  projective (inverse) limits.\\

\begin{bk}\label{I: 3.2.1} Let $X$ be a topological space, $\mathfrak{B}$ a basis for the opens of $X$ and $\emph{\textbf{K}}$ a category. Consider the (full) subcategory of $\ou(X)$ whose objects are the elements of $\mathfrak{B}$. Every contravariant functor $\mathscr{F}$ from this subcategory to $\emph{\textbf{K}}$ is called a presheaf over $\mathfrak{B}$ (or, a $\mathfrak{B}$-presheaf) with values in $\emph{\textbf{K}}$. Therefore there is a family of objects $\mathscr{F}(U)\in\emph{\textbf{K}}$, attached to each $U\in\mathfrak{B}$, and a family of morphisms $\rho_{U}^{V}:\mathscr{F}(V)\rightarrow\mathscr{F}(U)$ defined for each pair $(U,V)$ of elements of $\mathfrak{B}$ so that $U\subseteq V$, with the conditions that $\rho_{U}^{U}$ is the identity morphism and $\rho_{U}^{W}=\rho_{U}^{V}\circ\rho_{V}^{W}$ for any $U,V,W\in\mathfrak{B}$ with $U\subseteq V\subseteq W$.
Corresponding with such $\mathscr{F}$ one can associate a presheaf (in the usual sense) $\mathscr{F}'$ over $X$ with values in the category $\emph{\textbf{K}}$ as follows. For each open subspace $U$ of $X$, take $\mathscr{F}'(U)=\lim\limits_{V}\mathscr{F}(V)$ to be the projective limit of the family $(\mathscr{F}(V))_{V}$ where $V$ runs through the set of open subspaces belong to $\mathfrak{B}$ which are contained in $U$ (note that ordered by inclusion this set is not necessarily a filtered set; in fact if each such $V$ is nonempty then this set is filtered if and only if $U$ is either irreducible or empty), and for each two open subspaces $U$ and $U'$ of $X$ with $U\subseteq U'$, we define $(\rho')_{U}^{U'}:\mathscr{F}'(U')\rightarrow\mathscr{F}'(U)$ to be the unique morphism in $\emph{\textbf{K}}$ as the projective limit of the morphisms $\mathscr{F}'(U')\rightarrow\mathscr{F}(V)$ where $V$ runs through the set of open subspaces belong to $\mathfrak{B}$ which are contained in $U$. By the universal property of the projective limit, $(\rho')_{U}^{U}$ is the identity morphism and $(\rho')_{U}^{U''}=(\rho')_{U}^{U'}\circ(\rho')_{U'}^{U''}$ for any triple
$(U, U', U'')$ of open subspaces of $X$ with $U\subseteq U'\subseteq U''$. Moreover, for each basis open subspace $U\in\mathfrak{B}$, the canonical morphism $\can_{U}:\mathscr{F}'(U)\rightarrow\mathscr{F}(U)$ is an isomorphism, which permits us to identify these two objects. Because by the universal property of the projective limit there exists a (unique) morphism $\eta_{U}:\mathscr{F}(U)\rightarrow\mathscr{F}'(U)$ so that  for each $V\in\mathfrak{B}$ with $V\subseteq U$, $\rho_{V}^{U}=\can_{V}\circ\eta_{U}$ where $\can_{V}:\mathscr{F}'(U)\rightarrow\mathscr{F}(V)$ is the canonical morphism (sometimes we denote it also by $\can_{U, V}$ to avoid any possible confusion which may arise), in particular $\Id=\rho_{U}^{U}=\can_{U}\circ\eta_{U}$, also for every such open subspace $V$, $\can_{V}\circ(\eta_{U}\circ\can_{U})=\can_{V}$ and so $\eta_{U}\circ\can_{U}=\Id$, hence $\can_{U}$ is an isomorphism.\\
\end{bk}

\begin{bk}\label{I: 3.2.2} The presheaf $\mathscr{F}'$ as defined in the above is a sheaf if and only if the presheaf $\mathscr{F}$ over $\mathfrak{B}$ satisfies in the following condition:\\

($\textbf{F}_{0}$) For each covering $(U_{\alpha})$ of $U\in\mathfrak{B}$ by the open subspaces $U_{\alpha}\in\mathfrak{B}$ contained in $U$, and for each object $T\in\emph{\textbf{K}}$, the function which matches each $f\in\Hom_{\emph{\textbf{K}}}(T, \mathscr{F}(U))$ to the family $(\rho_{U_{\alpha}}^{U}\circ f)\in\prod\limits_{\alpha}\Hom_{\emph{\textbf{K}}}(T, \mathscr{F}(U_{\alpha}))$ is a bijective from $\Hom_{\emph{\textbf{K}}}(T, \mathscr{F}(U))$ onto the set of families $(f_{\alpha})\in\prod\limits_{\alpha}\Hom_{\emph{\textbf{K}}}(T, \mathscr{F}(U_{\alpha}))$ so that $\rho_{V}^{U_{\alpha}}\circ f_{\alpha}=\rho_{V}^{U_{\beta}}\circ f_{\beta}$ for each pair of indices $(\alpha, \beta)$ and for each $V\in\mathfrak{B}$ with $V\subseteq U_{\alpha}\cap U_{\beta}$.\\

Because first assume that $\mathscr{F}'$ is a sheaf. Take $\rho'_{\alpha}=(\rho')^{U}_{U_{\alpha}}: \mathscr{F}'(U)\rightarrow\mathscr{F}'(U_{\alpha})$ and $\rho'_{\alpha\beta}=(\rho')_{U_{\alpha}\cap U_{\beta}}^{U_{\alpha}}:\mathscr{F}'(U_{\alpha})
\rightarrow\mathscr{F}'(U_{\alpha}\cap U_{\beta})$, $\mathscr{F}$ satisfies in the condition ($\textbf{F}_{0}$), because if there exist two morphisms $f,g: T\rightarrow\mathscr{F}(U)$ so that for each $\alpha$, $\rho_{\alpha}\circ f=\rho_{\alpha}\circ g$, then $\rho'_{\alpha}\circ(\eta_{U}\circ f)=\rho'_{\alpha}\circ(\eta_{U}\circ g)$ since the following diagram is commutative $$\xymatrix{
\mathscr{F}(U) \ar[r]^{\eta_{U}} \ar[d]^{\rho_{\alpha}} & \mathscr{F}'(U) \ar[d]^{\rho'_{\alpha}} \\ \mathscr{F}(U_{\alpha})\ar[r]^{\eta_{U_{\alpha}}} & \mathscr{F}'(U_{\alpha})} $$ hence $\eta_{U}\circ f=\eta_{U}\circ g$ and so $f=g$ since $\eta_{U}$ is a monomorphism (in fact an isomorphism). Now let $(f_{\alpha})\in\prod\limits_{\alpha}\Hom_{\emph{\textbf{K}}}(T, \mathscr{F}(U_{\alpha}))$ be a family so that $\rho_{V}^{U_{\alpha}}\circ f_{\alpha}=\rho_{V}^{U_{\beta}}\circ f_{\beta}$ for each pair of indices $(\alpha, \beta)$ and for each $V\in\mathfrak{B}$ such that $V\subseteq U_{\alpha}\cap U_{\beta}$, then the family $(\eta_{U_{\alpha}}\circ f_{\alpha})\in\prod\limits_{\alpha}\Hom_{\emph{\textbf{K}}}(T, \mathscr{F}'(U_{\alpha}))$ satisfies in the condition $\rho'_{\alpha\beta}\circ(\eta_{U_{\alpha}}\circ f_{\alpha})=\rho'_{\beta\alpha}\circ(\eta_{U_{\beta}}\circ f_{\beta})$, because $(\rho')_{V}^{U_{\alpha}\cap U_{\beta}}\circ\rho'_{\alpha\beta}\circ(\eta_{U_{\alpha}}\circ f_{\alpha})=(\rho')_{V}^{U_{\alpha}\cap U_{\beta}}\circ\rho'_{\beta\alpha}\circ(\eta_{U_{\beta}}\circ f_{\beta})$
for every open subspace $V\in\mathfrak{B}$ contained in $U_{\alpha}\cap U_{\beta}$ (in particular for an open covering $(V_{\gamma})$ for
$U_{\alpha}\cap U_{\beta}$ formed by $V_{\gamma}\in\mathfrak{B}$ with $V_{\gamma}\subseteq U_{\alpha}\cap U_{\beta}$ ), and so there exists a morphism $f\in\Hom_{\emph{\textbf{K}}}(T, \mathscr{F}'(U))$ so that for each $\alpha$, $\rho'_{\alpha}\circ f=\eta_{U_{\alpha}}\circ f_{\alpha}$, then the morphism $\can_{U}\circ f\in\Hom_{\emph{\textbf{K}}}(T, \mathscr{F}(U))$ satisfies in the condition $\rho_{\alpha}\circ(\can_{U}\circ f)=f_{\alpha}$ for each $\alpha$.\\
Conversely, suppose that $\mathscr{F}$ satisfies in the condition ($\textbf{F}_{0}$). We separate the proof into two steps.\\
Step 1. Let $\mathfrak{B}'$ be a basis for the topology of $X$ which is contained in $\mathfrak{B}$ and denote by $\mathscr{F}_{0}$ the $\mathfrak{B}'-$presheaf induced by $\mathscr{F}$, i.e., for each open subspace $V\in\mathfrak{B}'$,  $\mathscr{F}_{0}(V)=\mathscr{F}(V)$ and for any $V, U\in\mathfrak{B}'$ with $V\subseteq U$, $(\rho_{0})_{V}^{U}=\rho_{V}^{U}$, denote by $\mathscr{F}''$  the presheaf (in the usual sense) corresponding to $\mathscr{F}_{0}$ as introduced in \ref{I: 3.2.1}; then $\mathscr{F}''$ is canonically isomorphic to $\mathscr{F}'$. Because, for each open subspace $U$ of $X$, take the morphism $\zeta_{U}:\mathscr{F}'(U)\rightarrow\mathscr{F}''(U)$ to be the projective limit of the canonical morphisms $\can_{V}: \mathscr{F}'(U)\rightarrow\mathscr{F}(V)=\mathscr{F}_{0}(V)$ where $V$ runs through the set of open subspaces belong to $\mathfrak{B}'$ contained in $U$; by the universal property of the projective limit $\zeta:\mathscr{F}'\rightarrow\mathscr{F}''$ is actually a morphism of presheaves. To construct $\xi:\mathscr{F}''\rightarrow\mathscr{F}'$, the inverse morphism of $\zeta$, we act as follows. For each fixed $V\in\mathfrak{B}$ with $V\subseteq U$, and for each $V_{\alpha}\in\mathfrak{B}'$ with $V_{\alpha}\subseteq V$, the family of morphisms $(\can'_{V_{\alpha}})$ where $\can'_{V_{\alpha}}: T=\mathscr{F}''(U)\rightarrow\mathscr{F}_{0}(V_{\alpha})
=\mathscr{F}(V_{\alpha})$ is the canonical morphism satisfies in the condition $\rho_{W}^{V_{\alpha}}\circ\can'_{V_{\alpha}}
=\rho_{W}^{V_{\beta}}\circ\can'_{V_{\beta}}$ for each pair of indices $(\alpha, \beta)$ and for each $W\in\mathfrak{B}$ with $W\subseteq V_{\alpha}\cap V_{\beta}$ (because for an open covering  $W=\bigcup\limits_{\gamma}W_{\gamma}$ with $W_{\gamma}\in\mathfrak{B}'$ for each $\gamma$, one has $\rho_{W_{\gamma}}^{W}\circ\rho_{W}^{W_{\alpha}}\circ
\can'_{V_{\alpha}}=\can'_{W_{\gamma}}=\rho_{W_{\gamma}}^{W}
\circ\rho_{W}^{W_{\beta}}\circ\can'_{V_{\beta}}$),
hence there exists a unique morphism $\omega_{V}: T=\mathscr{F}''(U)\rightarrow\mathscr{F}(V)$ so that $\can'_{V_{\alpha}}=\rho_{V_{\alpha}}^{V}\circ\omega_{V}$ and so $\omega_{V'}=\rho_{V'}^{V}\circ\omega_{V}$ for each pair $(V', V)$ of elements of $\mathfrak{B}$ with $V'\subseteq V$. Let $\xi_{U}:\mathscr{F}''(U)\rightarrow\mathscr{F}'(U)$ be the projective limit of the morphisms $\omega_{V}$ where $V$ runs through the set of open subspaces belong to $\mathfrak{B}$ which are contained in $U$. One can check that $\xi:\mathscr{F}''\rightarrow\mathscr{F}'$ is a morphism of presheaves and also it is the inverse of $\zeta$.\\
Step 2. Let $(U_{\alpha})$ be an open covering of an open subspace $U$ of $X$ each one is contained in $U$, let $\mathfrak{B}'$ be the subfamily of $\mathfrak{B}$ consisted by the open subspaces $V\in\mathfrak{B}$ which are contained at least in some $U_{\alpha}$, obviously $\mathfrak{B}'$ is a basis for the subspace topology of $U$ which is contained in the basis $\mathfrak{B}_{U}=\{V\in\mathfrak{B} : V\subseteq U\}$ of $U$, $\mathscr{F}'|_{U}$ is the presheaf corresponding to the presheaf $\mathscr{F}|_{U}$ over $\mathfrak{B}_{U}$, let $\mathscr{F}''$ be the presheaf (over $U$) corresponding to the presheaf $\mathscr{F}_{0}$ over $\mathfrak{B}'$ where $\mathscr{F}_{0}$ is induced via $\mathscr{F}|_{U}$. To verify the axiom  ($\textbf{F}$), first let that $f, g:T\rightarrow\mathscr{F}'(U)$ be two morphisms so that for each $\alpha$, $(\rho')_{U_{\alpha}}^{U}\circ f=(\rho')_{U_{\alpha}}^{U}\circ g$, then for each $V\in\mathfrak{B}'$ there exists some $\beta$ so that $V\subseteq U_{\beta}$, and we have $\can'_{V}\circ\zeta_{U}\circ f=\can_{U, V}\circ f=\can_{U_{\beta},V}\circ(\rho')_{U_{\beta}}^{U}\circ f=\can'_{V}\circ\zeta_{U}\circ g$ so $\zeta_{U}\circ f=\zeta_{U}\circ g$ this implies that $f=g$.  Now let $(f_{\alpha})\in\prod\limits_{\alpha}\Hom_{\emph{\textbf{K}}}(T, \mathscr{F}'(U_{\alpha}))$ be a family so that for each pair of indices $(\alpha, \beta)$, $\rho'_{\alpha\beta}\circ f_{\alpha}=\rho'_{\beta\alpha}\circ f_{\beta}$. For each $V\in\mathfrak{B}'$ there exists some $\beta$ so that $V\subseteq U_{\beta}$, take $\theta_{V}=\can_{U_{\beta}, V}\circ f_{\beta}$ which is independent of choosing such $\beta$, then for each pair $(V', V)$ of elements of $\mathfrak{B}'$ with $V'\subseteq V$,  $\rho_{V'}^{V}\circ\theta_{V}=\theta_{V'}$, hence there exist a unique morphism $\theta: T\rightarrow\mathscr{F}''(U)$ so that for each $V\in\mathfrak{B}'$, $\can'_{U, V}\circ\theta=\theta_{V}$, then the morphism $\xi_{U}\circ\theta:T\rightarrow\mathscr{F}'(U)$ satisfies in the condition $f_{\alpha}=(\rho')_{U_{\alpha}}^{U}\circ(\xi_{U}\circ\theta)$ for each $\alpha$, because as we showed above, the map $\Hom_{\emph{\textbf{K}}}(T, \mathscr{F}'(U_{\alpha}))\rightarrow
\prod\limits_{\gamma}\Hom_{\emph{\textbf{K}}}(T, \mathscr{F}'(V_{\alpha\gamma}))$ given by $f\rightsquigarrow((\rho')_{V_{\alpha\gamma}}^{U_{\alpha}})\circ f)$ is injective where the family $(V_{\alpha\gamma})_{\gamma}$
is an open covering for $U_{\alpha}$ formed by elements of $\mathfrak{B}'$.\\
The $\mathfrak{B}-$presheaf  $\mathscr{F}$, by abuse of the language, is said to be a $\mathfrak{B}-$sheaf or a sheaf over $\mathfrak{B}$ if $\mathscr{F}$ satisfies in the axiom $(\textbf{F}_{0})$.\\
\end{bk}

\begin{bk}\label{I: 3.2.3} Let $\mathscr{F}, \mathscr{G}$ be two presheaves over the basis $\mathfrak{B}$ with values in the category $\emph{\textbf{K}}$, we define a morphism $u:\mathscr{F}\rightarrow\mathscr{G}$ as a family $(u_{V})_{V\in\mathfrak{B}}$ of morphisms $u_{V}:\mathscr{F}(V)\rightarrow\mathscr{G}(V)$ so that for each pair $(W, V)$ of elements of $\mathfrak{B}$ with $W\subseteq V$, the following diagram is commutative  $$\xymatrix{
\mathscr{F}(V) \ar[r]^{u_{V}} \ar[d]^{\rho_{W}^{V}} & \mathscr{G}(V) \ar[d]^{\varrho_{W}^{V}} \\ \mathscr{F}(W)\ar[r]^{u_{W}} & \mathscr{G}(W) .} $$ Taking into account the notations of \ref{I: 3.2.1}, corresponding with such morphism $u$, one can then deduce the morphism $u':\mathscr{F}'\rightarrow\mathscr{G}'$ of presheaves (ordinary) by taking $u'_{U}$, for each open subspace $U$ of $X$, to be the unique morphism in $\emph{\textbf{K}}$ as the projective limit of the morphisms $u_{V}\circ\can_{V}:\mathscr{F}'(U)\rightarrow\mathscr{G}(V)$ where $V$ runs through the set of elements of $\mathfrak{B}$ contained in $U$. For each pair $(U', U)$ of open subspaces of $X$ with $U'\subseteq U$, by the universal property of the projective limit, the following diagram is commutative $$\xymatrix{
\mathscr{F}'(U) \ar[r]^{u'_{U}} \ar[d]^{(\rho')_{U'}^{U}} & \mathscr{G}'(U) \ar[d]^{(\varrho')_{U'}^{U}} \\ \mathscr{F}'(U')\ar[r]^{u'_{U'}} & \mathscr{G}'(U')}$$ hence $u'$ is actually a morphism of presheaves.\\
\end{bk}

\begin{bk}\label{I: 3.2.4} Let $\mathscr{F}$ be a presheaf over the basis $\mathfrak{B}$ with values in the category $\emph{\textbf{K}}$, for each point $x\in X$, the set of neighborhoods of $x$ belong to $\mathfrak{B}$ is a filtered (directed) subset of the set of open neighborhoods of $x$, if $\mathscr{F}'$ is the presheaf (ordinary) corresponding to $\mathscr{F}$ then the fiber $(\mathscr{F}')_{x}$ is canonically isomorphic to the $\lim\limits_{x\in V}\mathscr{F}(V)$ where $V$ runs through the set of elements of $\mathfrak{B}$ containing $x$ (because, more generally, let $(A_{i}, \phi_{i, j})_{i,j\in I}$ be an inductive (direct) system formed by objects and morphisms of $\emph{\textbf{K}}$ indexed by the filtered set $I$, if $J$ is a filtered subset of $I$ so that for each $i\in I$ there exists some $j\in J$ in which $i\leq j$, then the object $\lim\limits_{i\in I}A_{i}$ is canonically isomorphic to the object  $\lim\limits_{j\in J}A_{j}$).
\end{bk}

\begin{bk}\label{I: 3.2.5} Let $\mathscr{F}$ be a sheaf (ordinary) with values in $\emph{\textbf{K}}$, $\mathscr{F}_{1}$ the sheaf over $\mathfrak{B}$ obtained by the restriction of $\mathscr{F}$ over $\mathfrak{B}$, then the sheaf  $\mathscr{F}'_{1}$ (corresponding with $\mathscr{F}_{1}$), by virtue of the condition $(\textbf{F})$ and the universal property of the projective limit, is canonically isomorphic to $\mathscr{F}$. Because, for each fixed open subspace $U$ of $X$, there exists a unique morphism $\theta_{U}:\mathscr{F}(U)\rightarrow\mathscr{F}_{1}'(U)$ in which  $\rho_{V}^{U}=\can_{U, V}\circ\theta_{U}$ for each open subspace $V\in\mathfrak{B}$ contained in $U$, using then the universal property of the projective limit, $\theta:\mathscr{F}\rightarrow\mathscr{F}_{1}'$ is actually a morphism of presheaves. On the other hand, the family of morphisms $(\can_{U, V})\in\prod\limits_{V}
\Hom_{\emph{\textbf{K}}}(\mathscr{F}_{1}'(U), \mathscr{F}_{1}(V))$ where $V$ runs through the set of elements of $\mathfrak{B}$ contained in $U$, satisfies in the condition $\rho_{V\cap V'}^{V}\circ\can_{U,V}=\rho_{V\cap V'}^{V'}\circ\can_{U,V'}$ for each pair $(V, V')$ of elements of $\mathfrak{B}$ with $V, V'\subseteq U$ (because for an open covering $V\cap V'=\bigcup\limits_{W\in\mathfrak{B}}W$, one has $\rho_{W}^{V\cap V'}\circ\rho_{V\cap V'}^{V}\circ\can_{U,V}=\can_{U,W}=\rho_{W}^{V\cap V'}\circ\rho_{V\cap V'}^{V'}\circ\can_{U,V'}$, then the desired assertion implies from the condition $(\textbf{F})$ ), hence by the condition $(\textbf{F})$ there exists a (unique) morphism $\psi_{U}:\mathscr{F}_{1}'(U)\rightarrow\mathscr{F}(U)$ so that $\can_{U,V}=\rho_{V}^{U}\circ\psi_{U}$ for each $V\in\mathfrak{B}$ contained in $U$,  using the condition $(\textbf{F})$ once more we then conclude that $\psi:\mathscr{F}_{1}'\rightarrow\mathscr{F}$ is actually a morphism of presheaves. Finally, it is easy to observe that the morphisms $\psi\circ\theta$ and $\theta\circ\psi$ are indeed the identity morphisms.\\
Let $\mathscr{G}$ be a second sheaf (ordinary) over $X$ with values in $\emph{\textbf{K}}$, and $u:\mathscr{F}\rightarrow\mathscr{G}$ a morphism, the previous remark (\ref{I: 3.2.3}) shows that given the morphisms $u_{V}:\mathscr{F}(V)\rightarrow\mathscr{G}(V)$ for each $V\in\mathfrak{B}$, uniquely determines $u$ completely (i.e. the morphism $u$ can be recovered with having in hand the family of morphisms $(u_{V})_{V\in\mathfrak{B}}$). \\
\end{bk}

\begin{bk}\label{I: 3.2.6}  The category of sheaves over $X$ with values in $\emph{\textbf{K}}$, denoted by $\sh_{\emph{\textbf{K}}}(X)$, admits the projective limits. Because suppose that $(\mathscr{F}_{\lambda}, u_{\lambda\mu})_{\lambda,\mu\in\Lambda}$ is a projective system of sheaves over $X$ with values in $\emph{\textbf{K}}$ indexed by a partially ordered set $\Lambda$ (which is not necessarily filtered), to construct the projective limit of this system we act as follows. For each fixed open subspace $U$ of $X$, the family $(\mathscr{F}_{\lambda}(U), (u_{\lambda\mu})_{U})$ is a projective system over the indexed set $\Lambda$
formed by objects and morphisms of $\emph{\textbf{K}}$, take $\mathscr{F}(U)=
\lim\limits_{\lambda\in\Lambda}\mathscr{F}_{\lambda}(U)$ to be the projective limit of this system, and for each pair $(U, U')$ of open subsets of $X$ with $U'\subseteq U$, take $\rho_{U'}^{U}:\mathscr{F}(U)\rightarrow\mathscr{F}(U')$ to be the projective limit of the morphisms $(\rho_{\lambda})_{U'}^{U}\circ\can_{U,\lambda}:
\mathscr{F}(U)\rightarrow\mathscr{F}_{\lambda}(U')$ where $\lambda$ runs through the indexed set $\Lambda$ and $\can_{U,\lambda}:\mathscr{F}(U)\rightarrow\mathscr{F}_{\lambda}(U)$ is the canonical morphism. Obviously $\mathscr{F}$ is a presheaf over $X$ with values in $\emph{\textbf{K}}$. Let $(U_{\alpha})$ be an open covering of $U$ formed by open subspaces contained in $U$ and let $T$ be an arbitrary object of $\emph{\textbf{K}}$, to verify the condition $(\textbf{F})$ first suppose that there exist two morphisms $f,g:T\rightarrow\mathscr{F}(U)$ so that $\rho_{\alpha}\circ f=\rho_{\alpha}\circ g$ for each $\alpha$. For each $\lambda\in\Lambda$, set $\theta_{\lambda}=\can_{U,\lambda}\circ f: T\rightarrow\mathscr{F}_{\lambda}(U)$, then for each pair $(\lambda, \mu)$ of elements of $\Lambda$ with $\lambda\leq\mu$, one has $\theta_{\mu}=(u_{\lambda\mu})_{U}\circ\theta_{\lambda}:
T\rightarrow\mathscr{F}_{\mu}(U)$ (because $\mathscr{F}_{\mu}$ is a sheaf and hence satisfies in the condition  $(\textbf{F})$), hence by the universal property of the projective limit there exists a unique morphism $\theta:T\rightarrow\mathscr{F}(U)$ so that for each $\lambda$, $\theta_{\lambda}=\can_{U,\lambda}\circ\theta$, but both of the morphisms $f$ and $g$ do the task of $\theta$ therefore $f=g$. Now let $(f_{\alpha})\in\prod\limits_{\alpha}\Hom_{\emph{\textbf{K}}}(T, \mathscr{F}(U_{\alpha}))$ be a family of morphisms so that $\rho_{\alpha\beta}\circ f_{\alpha}=\rho_{\beta\alpha}\circ f_{\beta}$ for each pair of indices $(\alpha, \beta)$; for each $\lambda\in\Lambda$ and for each $\alpha$ set $h_{\lambda,\alpha}=\can_{U_{\alpha},\lambda}\circ f_{\alpha}:T\rightarrow\mathscr{F}_{\lambda}(U_{\alpha})$, then for each fixed $\lambda$, the family of morphisms $(h_{\lambda,\alpha})_{\alpha}$ satisfies in the condition $(\rho_{\lambda})_{\alpha\beta}\circ h_{\lambda,\alpha}=(\rho_{\lambda})_{\beta\alpha}\circ h_{\lambda,\beta}$ for each pair of indices $(\alpha,\beta)$, hence there exists a unique morphism $h_{\lambda}:T\rightarrow\mathscr{F}_{\lambda}(U)$ so that for each $\alpha$, $h_{\lambda,\alpha}=(\rho_{\lambda})_{\alpha}\circ h_{\lambda}$. For each pair of indices $(\lambda, \mu)$ of elements of $\Lambda$ with $\lambda\leq\mu$, we have $h_{\mu}=(u_{\lambda\mu})_{U}\circ h_{\lambda}$, hence by the universal property of the projective limit there exists a unique morphism $h:T\rightarrow\mathscr{F}(U)$ so that for each $\lambda$, $h_{\lambda}=\can_{U,\lambda}\circ h$; and eventually for each $\alpha$ we have $f_{\alpha}=\rho_{\alpha}\circ h$, hence $\mathscr{F}$ is actually a sheaf over $X$ with values in $\emph{\textbf{K}}$.\\
For each $\lambda$, $\can_{-,\lambda}:\mathscr{F}\rightarrow\mathscr{F}_{\lambda}$ is actually a morphism of sheaves and for each pair $(\lambda,\mu)$ of elements of $\Lambda$ with $\lambda\leq\mu$, $\can_{-,\mu}=u_{\lambda\mu}\circ\can_{-,\lambda}$; now let $\mathscr{G}$ be an arbitrary sheaf over $X$ with values in $\emph{\textbf{K}}$ together with the family $(\omega_{\lambda}:\mathscr{G}
\rightarrow\mathscr{F}_{\lambda})$ of morphisms (of sheaves) in which for each pair $(\lambda,\mu)$ of elements of $\Lambda$ with $\lambda\leq\mu$, $\omega_{\mu}=u_{\lambda\mu}\circ\omega_{\lambda}$; for each open subspace $U$ of $X$, by the universal property of the projective limit, there exists a unique morphism $\omega_{U}:\mathscr{G}(U)\rightarrow\mathscr{F}(U)$ so that for each $\lambda$, $(\omega_{\lambda})_{U}=\can_{U,\lambda}\circ\omega_{U}$; obviously $\omega:\mathscr{G}\rightarrow\mathscr{F}$ is a morphism of sheaves and satisfies in the condition $\omega_{\lambda}=\can_{-,\lambda}\circ\omega$ for each $\lambda$, hence the pair $(\mathscr{F},(\can_{-,\lambda})_{\lambda})$ is the projective limit of the system  $(\mathscr{F}_{\lambda}, u_{\lambda\mu})_{\lambda,\mu\in\Lambda}$.\\
\end{bk}

\section{Gluing of sheaves}

Throughout of this section, the category $\emph{\textbf{K}}$ admits the projective limits.\\

\begin{bk}\label{I: 3.3.1} Let $X$ be a topological space, $\mathscr{U}=(U_{\lambda})_{\lambda\in L}$ an open covering for $X$, and for each $\lambda\in L$, $\mathscr{F}_{\lambda}$ is a sheaf over the subspace $U_{\lambda}$ with values in $\emph{\textbf{K}}$; for each couple of indices $(\lambda,\mu)$, suppose that an isomorphism $\theta_{\lambda,\mu}:\mathscr{F}_{\mu}|_{U_{\lambda}\cap U_{\mu}}\rightarrow\mathscr{F}_{\lambda}|_{U_{\lambda}\cap U_{\mu}}$ is given in which $\theta_{\lambda,\lambda}$ is the identity morphism and for each triple of indices $(\lambda,\mu,\nu)$, we have $\theta'_{\lambda,\nu}=\theta'_{\lambda,\mu}\circ\theta'_{\mu,\nu}$  (gluing conditions of the family $(\theta_{\lambda,\mu})_{\lambda,\mu\in L}$ ) where the restriction of $\theta_{\lambda,\mu}$ over $U_{\lambda}\cap U_{\mu}\cap U_{\nu}$ is denoted by $\theta'_{\lambda,\mu}$. Then there exists a sheaf $\mathscr{F}'$ over $X$ with values in $\emph{\textbf{K}}$ and for each $\lambda$, an isomorphism $\eta_{\lambda}:\mathscr{F}'|_{U_{\lambda}}
\rightarrow\mathscr{F}_{\lambda}$ such that for each couple of indices $(\lambda,\mu)$, we have $\theta_{\lambda,\mu}=\eta'_{\lambda}\circ(\eta'_{\mu})^{-1}$ where the restriction of $\eta_{\lambda}$ over $U_{\lambda}\cap U_{\mu}$ is denoted by $\eta'_{\lambda}$. Moreover, under these conditions, the sheaf $\mathscr{F}'$ and the family $(\eta_{\lambda})_{\lambda\in L}$, with approximation of a unique isomorphism, are unique. Because let $\mathfrak{B}$ be the set of all open subspaces of $X$ which are contained at least in some $U_{\lambda}$, obviously $\mathfrak{B}$ is a basis for the topology of $X$; for each $V\in\mathfrak{B}$, denote by $\mathscr{A}_{V}$ the set of all $\lambda\in L$ in which $V\subseteq U_{\lambda}$, which is a nonempty subset; by the axiom of choice we can take a choice function $\tau:\mathfrak{B}\rightarrow
\bigcup\limits_{V\in\mathfrak{B}}\mathscr{A}_{V}$; we then define $\mathscr{F}(V):=\mathscr{F}_{\tau V}(V)$ and also for each pair $(V', V)$ of elements of $\mathfrak{B}$ with $V'\subseteq V$, we define $\rho_{V'}^{V}:=(\theta_{\tau V',\tau V})_{V'}\circ(\rho_{\tau V})_{V'}^{V}:\mathscr{F}(V)\rightarrow\mathscr{F}(V')$, obviously $\mathscr{F}$ is a presheaf over $\mathfrak{B}$ with values in $\emph{\textbf{K}}$ (because for any triple $(V,V',V'')$ of elements of $\mathfrak{B}$ with $V''\subseteq V'\subseteq V$, we have $\rho_{V''}^{V'}\circ\rho_{V'}^{V}=((\theta_{\tau V'',\tau V'})_{V''}\circ(\rho_{\tau V'})_{V''}^{V'})\circ((\theta_{\tau V',\tau V})_{V'}\circ(\rho_{\tau V})_{V'}^{V})=
(\theta_{\tau V'',\tau V'})_{V''}\circ(\theta_{\tau V',\tau V})_{V''}\circ(\rho_{\tau V})_{V''}^{V'}\circ(\rho_{\tau V})_{V'}^{V}=(\theta_{\tau V'',\tau V})_{V''}\circ(\rho_{\tau V})_{V''}^{V}=\rho_{V''}^{V}$).  Let $(V_{\alpha})$ be an open covering of $V\in\mathfrak{B}$ formed by open subspaces belong to $\mathfrak{B}$ which are contained in $V$, and let $T$ be an arbitrary object of $\emph{\textbf{K}}$, we verify the condition $(\textbf{F}_{0})$ for the presheaf  $\mathscr{F}$; first suppose that there exist two morphisms $f,g\in\Hom_{\emph{\textbf{K}}}(T,\mathscr{F}(V))$ so that for each $\alpha$, $\rho_{\alpha}\circ f=\rho_{\alpha}\circ g$, this implies that $(\rho_{\tau V})_{V_{\alpha}}^{V}\circ f=(\rho_{\tau V})_{V_{\alpha}}^{V}\circ g$ since $(\theta_{\tau V_{\alpha},\tau V})_{V_{\alpha}}$ is a monomorphism (in fact it is an isomorphism), but $\mathscr{F}_{\tau V}$ is a sheaf over $U_{\tau V}$ hence $f=g$. Now let $(f_{\alpha})\in\prod\limits_{\alpha}
\Hom_{\emph{\textbf{K}}}(T,\mathscr{F}(V_{\alpha}))$ be a family of morphisms in which  for each pair of indices $(\alpha,\beta)$ and for each $W\in\mathfrak{B}$ with $W\subseteq V_{\alpha}\cap V_{_{\beta}}$, we have
$\rho_{W}^{V_{\alpha}}\circ f_{\alpha}=\rho_{W}^{V_{\beta}}\circ f_{\beta}$. For each $\alpha$, set $g_{\alpha}=(\theta_{\tau V_{\alpha},\tau V})_{V_{\alpha}}^{-1}\circ f_{\alpha}:T\rightarrow\mathscr{F}_{\tau V}(V_{\alpha})$, then for each pair of indices $(\alpha,\beta)$, we have $(\rho_{\tau V})_{\alpha\beta}\circ g_{\alpha}=(\rho_{\tau V})_{\beta\alpha}\circ g_{\beta}$ (because for each $W\in\mathfrak{B}$ with $W\subseteq V_{\alpha}\cap V_{_{\beta}}$, we have $(\rho_{\tau V})_{W}^{V_{\alpha}\cap V_{\beta}}\circ(\rho_{\tau V})_{\alpha\beta}\circ g_{\alpha}=(\rho_{\tau V})_{W}^{V_{\alpha}}\circ g_{\alpha}=(\theta_{\tau V_{\alpha},\tau V})_{W}^{-1}\circ(\rho_{\tau V_{\alpha}})_{W}^{V_{\alpha}}\circ f_{\alpha}=(\theta_{\tau W,\tau V})_{W}^{-1}\circ(\theta_{\tau W,\tau V_{\alpha}})_{W}\circ(\rho_{\tau V_{\alpha}})_{W}^{V_{\alpha}}\circ f_{\alpha}=(\theta_{\tau W,\tau V})_{W}^{-1}\circ(\theta_{\tau W,\tau V_{\beta}})_{W}\circ(\rho_{\tau V_{\beta}})_{W}^{V_{\beta}}\circ f_{\beta}=(\theta_{\tau V_{\beta},\tau V})_{W}^{-1}\circ(\rho_{\tau V_{\beta}})_{W}^{V_{\beta}}\circ f_{\beta}=(\rho_{\tau V})_{W}^{V_{\beta}}\circ(\theta_{\tau V_{\beta},\tau V})_{V_{\beta}}^{-1}\circ f_{\beta}=(\rho_{\tau V})_{W}^{V_{\beta}}\circ g_{\beta}=(\rho_{\tau V})_{W}^{V_{\alpha}\cap V_{\beta}}\circ(\rho_{\tau V})_{\beta\alpha}\circ g_{\beta}$; and $\mathscr{F}_{\tau V}$ is a sheaf over $U_{\tau V}$ ( instead of using this hypothesis, one can also take $W=V_{\alpha}\cap V_{\beta}\in\mathfrak{B}$)). Since  $\mathscr{F}_{\tau V}$ is a sheaf over $U_{\tau V}$ hence there exists a unique morphism $g:T\rightarrow\mathscr{F}_{\tau V}(V)=\mathscr{F}(V)$ so that for each $\alpha$, $g_{\alpha}=(\rho_{\tau V})_{\alpha}\circ g$; from this we observe that $\rho_{\alpha}\circ g=f_{\alpha}$ for each $\alpha$. Hence $\mathscr{F}$ is a sheaf over $\mathfrak{B}$, therefore the presheaf (ordinary) $\mathscr{F}'$ (corresponding with $\mathscr{F}$) is actually a sheaf over $X$ with values in $\emph{\textbf{K}}$. For each $\lambda\in L$ and for each open subspace $V$ of $U_{\lambda}$ (by the definition $V\in\mathfrak{B}$), the canonical morphism $\can_{V}:\mathscr{F}'(V)\rightarrow\mathscr{F}(V)$ is an isomorphism, take $(\eta_{\lambda})_{V}=(\theta_{\lambda,\tau V})_{V}\circ\can_{V}:\mathscr{F}'(V)\rightarrow\mathscr{F}_{\lambda}(V)$ then $\eta_{\lambda}:\mathscr{F}'|_{U_{\lambda}}\rightarrow\mathscr{F}_{\lambda}$ is actually an isomorphism of sheaves (because for any two open subspaces $V',V$ of $U_{\lambda}$ with $V'\subseteq V$, we have $(\rho_{\lambda})_{V'}^{V}\circ(\eta_{\lambda})_{V}=(\theta_{\lambda,\tau V})_{V'}\circ(\rho_{\tau V})_{V'}^{V}\circ\can_{V}=(\theta_{\lambda,\tau V'})_{V'}\circ(\theta_{\tau V',\tau V})_{V'}\circ(\rho_{\tau V})_{V'}^{V}\circ\can_{V}=(\theta_{\lambda,\tau V'})_{V'}\circ\rho_{V'}^{V}\circ\can_{V}=
(\theta_{\lambda,\tau V'})_{V'}\circ\can_{V,V'}=(\theta_{\lambda,\tau V'})_{V'}\circ\can_{V'}\circ(\rho')_{V'}^{V}=
(\eta_{\lambda})_{V'}\circ(\rho')_{V'}^{V}$ ), moreover for each couple of indices $(\lambda,\mu)$ of elements of $L$, and for each open subspace $V$ of $U_{\lambda}\cap U_{\mu}$, we have $(\eta_{\lambda}\circ\eta_{\mu}^{-1})_{V}=
(\eta_{\lambda})_{V}\circ(\eta_{\mu}^{-1})_{V}=(\theta_{\lambda,\mu})_{V}$, therefore $\theta_{\lambda,\mu}=\eta'_{\lambda}\circ(\eta'_{\mu})^{-1}$. \\ Finally, for uniqueness of the pair $(\mathscr{F}', (\eta_{\lambda})_{\lambda\in L})$, suppose that there is another pair $(\mathscr{G}, (\zeta_{\lambda})_{\lambda\in L})$ where $\mathscr{G}$ is a sheaf over $X$ with values in $\emph{\textbf{K}}$ and for each $\lambda\in L$, $\zeta_{\lambda}:
\mathscr{G}|_{U_{\lambda}}\rightarrow\mathscr{F}_{\lambda}$ is an isomorphism (of sheaves) in which for each couple of indices $(\lambda,\mu)$, $\theta_{\lambda,\mu}=\zeta'_{\lambda}\circ(\zeta'_{\mu})^{-1}$ where the restriction of $\zeta_{\lambda}$ over $U_{\lambda}\cap U_{\mu}$ is denoted by $\zeta'_{\lambda}$. (As another version of \ref{I: 3.2.5}, let $\mathscr{G}$ be a sheaf (ordinary) over $X$ with values in $\emph{\textbf{K}}$, $\mathscr{G}_{1}$ the restriction of $\mathscr{G}$ over $\mathfrak{B}$, $\mathscr{F}$ another sheaf over  $\mathfrak{B}$ with values in $\emph{\textbf{K}}$ so that there exists an isomorphism $\phi:\mathscr{G}_{1}\rightarrow\mathscr{F}$ as sheaves over  $\mathfrak{B}$, then there exists a unique isomorphism (of ordinary sheaves) $\Phi:\mathscr{G}\rightarrow\mathscr{F}'$ in which for each open subspace $U$ and for each $V\in\mathfrak{B}$ with $V\subseteq U$, we have  $\can_{U,V}\circ\Phi_{U}=\phi_{V}\circ\omega_{V}^{U}$ where $\mathscr{F}'$ is the sheaf corresponding with $\mathscr{F}$, $\omega_{V}^{U}:\mathscr{G}(U)\rightarrow\mathscr{G}(V)$ is the restriction morphism and $\can_{U,V}:\mathscr{F}'(U)\rightarrow\mathscr{F}(V)$ is the canonical morphism), hence as above, denote by $\mathscr{G}_{1}$ the restriction of $\mathscr{G}$ over $\mathfrak{B}$, moreover for each $V\in\mathfrak{B}$, define $\phi_{V}=(\zeta_{\tau V})_{V}:\mathscr{G}_{1}(V)=
\mathscr{G}(V)\rightarrow\mathscr{F}_{\tau V}(V)=\mathscr{F}(V)$, then
 $\mathscr{G}_{1}$ and $\mathscr{F}$ are isomorphic as sheaves over $\mathfrak{B}$ via $\phi:\mathscr{G}_{1}\rightarrow\mathscr{F}$, therefore there exists an (unique) isomorphism $\Phi:\mathscr{G}\rightarrow\mathscr{F}'$ (of ordinary sheaves) in which for each open subspace $U$ and for each $V\in\mathfrak{B}$ with $V\subseteq U$, we have  $\can_{U,V}\circ\Phi_{U}=\phi_{V}\circ\omega_{V}^{U}$, this in particular implies that for each $\lambda\in L$, $\zeta_{\lambda}=\eta_{\lambda}\circ\Phi_{U_{\lambda}}$; if there exists another isomorphism $\Psi:\mathscr{G}\rightarrow\mathscr{F}'$ having the property $\zeta_{\lambda}=\eta_{\lambda}\circ\Psi_{U_{\lambda}}$ for each $\lambda$, then for each $V\in\mathfrak{B}$, $\Phi_{V}=\Psi_{V}$, therefore $\Phi_{U}=\Psi_{U}$ for an arbitrary open subspace $U$ of $X$. \\
We say that the sheaf $\mathscr{F}'$ is obtained by gluing of the $\mathscr{F}_{\lambda}$ by means of the $\theta_{\lambda,\mu}$; the sheaves $\mathscr{F}'|_{U_{\lambda}}$ and $\mathscr{F}_{\lambda}$ are usually identified by means of $\eta_{\lambda}$.\\
Conversely, any sheaf $\mathscr{G}$ over $X$ with values in  $\emph{\textbf{K}}$ can be considered as the gluing of the sheaves $\mathscr{G}_{\lambda}=\mathscr{G}|_{U_{\lambda}}$ by means of the identity morphisms $\theta_{\lambda,\mu}$ (where $(U_{\lambda})$ is an open covering for $X$).\\
\end{bk}

\begin{bk}\label{I: 3.3.2} With the same notations as in \ref{I: 3.3.1}, for each $\lambda\in L$, let $\mathscr{G}_{\lambda}$ be a second sheaf over $U_{\lambda}$ with values in $\emph{\textbf{K}}$, and for each couple of indices $(\lambda,\mu)$ suppose that an isomorphism $\omega_{\lambda,\mu}:
\mathscr{G}_{\mu}|_{U_{\lambda}\cap U_{\mu}}\rightarrow\mathscr{G}_{\lambda}|_{U_{\lambda}\cap U_{\mu}}$ is given in which these isomorphisms satisfy in the gluing conditions. Suppose that the sheaf $\mathscr{G}'$ (together with the isomorphisms $\zeta_{\lambda}:
\mathscr{G}'|_{U_{\lambda}}\rightarrow\mathscr{G}_{\lambda}$ for each $\lambda$) is obtained by gluing of the $\mathscr{G}_{\lambda}$ by means of the $\omega_{\lambda,\mu}$, then there exists a functorial bijective given by $u\rightsquigarrow(\zeta_{\lambda}\circ u|_{U_{\lambda}}\circ\eta^{-1}_{\lambda})$
between the set $\Hom(\mathscr{F}',\mathscr{G}')$ and the set of families $(u_{\lambda})\in\prod\limits_{\lambda}
\Hom(\mathscr{F}_{\lambda},\mathscr{G}_{\lambda})$ so that for each pair of indices $(\lambda,\mu)$, the following diagram is commutative $$\xymatrix{
\mathscr{F}_{\mu}|_{U_{\lambda}\cap U_{\mu}} \ar[r]^{u_{\mu}} \ar[d]^{\theta_{\lambda,\mu}} & \mathscr{G}_{\mu}|_{U_{\lambda}\cap U_{\mu}} \ar[d]^{\omega_{\lambda,\mu}} \\ \mathscr{F}_{\lambda}|_{U_{\lambda}\cap U_{\mu}}\ar[r]^{u_{\lambda}} & \mathscr{G}_{\lambda}|_{U_{\lambda}\cap U_{\mu}} .} $$
Obviously the foregoing map is injective; now let $(u_{\lambda})\in\prod\limits_{\lambda}
\Hom(\mathscr{F}_{\lambda},\mathscr{G}_{\lambda})$ be a family so that for each pair of indices $(\lambda,\mu)$, the above diagram is commutative, suppose that $\sigma:\mathfrak{B}\rightarrow
\bigcup\limits_{V\in\mathfrak{B}}\mathscr{A}_{V}$ is a choice function for which $\mathscr{G}(V)=\mathscr{G}_{\sigma V}(V)$ for each $V\in\mathfrak{B}$, let $U$ be an arbitrary open subspace of $X$, for each $V\in\mathfrak{B}$ with $V\subseteq U$, take $\epsilon_{U,V}=(\omega_{\sigma V,\tau V})_{V}\circ (u_{\tau V})_{V}\circ\can_{V}\circ(\rho')_{V}^{U}:\mathscr{F}'(U)\rightarrow
\mathscr{G}(V)$, then one can observe that for each $V',V\in\mathfrak{B}$ with $V'\subseteq V\subseteq U$, $\epsilon_{U,V'}=\varrho_{V'}^{V}\circ\epsilon_{U,V}$ where $\varrho_{V'}^{V}:\mathscr{G}(V)\rightarrow\mathscr{G}(V')$ is the restriction morphism, hence by the universal property of the projective limit there exists a unique morphism $u_{U}:\mathscr{F}'(U)\rightarrow
\mathscr{G}'(U)$ in which for each $V\in\mathfrak{B}$ with $V\subseteq U$, $\epsilon_{U,V}=\can'_{U,V}\circ u_{U}$ where $\can'_{U,V}:\mathscr{G}'(U)\rightarrow\mathscr{G}(V)$ is the canonical morphism, $u:\mathscr{F}'\rightarrow
\mathscr{G}'$ is actually a morphism of sheaves and for each $\lambda$, $u_{\lambda}=\zeta_{\lambda}\circ u|_{U_{\lambda}}\circ\eta^{-1}_{\lambda}$.\\
\end{bk}

\begin{bk}\label{I: 3.3.3} With the same notations as in \ref{I: 3.3.1}, let $V$ be an open subspace of $X$, the sheaf $\mathscr{F}'|_{V}$ is the gluing of the $\mathscr{F}_{\lambda}|_{V\cap U_{\lambda}}$ by means of the $\theta_{\lambda,\mu}|_{V\cap U_{\lambda}\cap U_{\mu}}$.\\
\end{bk}

\section{Direct image of presheaves}

\begin{bk}\label{I: 3.4.1} Let $X$ and $Y$ be two topological spaces, $\psi:X\rightarrow Y$ a continuous map, $\mathscr{F}$ a presheaf over $X$ with values in a category $\emph{\textbf{K}}$, for each open subspace $U$ of $Y$ take $\psi_{\ast}(\mathscr{F})(U)=\mathscr{F}(\psi^{-1}(U))$, and for each two open subspaces $U, V$ of $Y$ with $U\subseteq V$, set $\omega_{U}^{V}=\rho_{\psi^{-1}(U)}^{\psi^{-1}(V)}:\psi_{\ast}
(\mathscr{F})(V)\rightarrow\psi_{\ast}(\mathscr{F})(U)$; obviously  $\psi_{\ast}(\mathscr{F})$ is a presheaf over $Y$ with values in $\emph{\textbf{K}}$, the presheaf $\psi_{\ast}(\mathscr{F})$ is called the direct image of  $\mathscr{F}$ by $\psi$; if $\mathscr{F}$ is a sheaf then obviously $\psi_{\ast}(\mathscr{F})$ is also a sheaf (over $Y$ with values in $\emph{\textbf{K}}$).\\
\end{bk}

\begin{bk}\label{I: 3.4.2} Let $\mathscr{F}_{1}$ and $\mathscr{F}_{2}$ be two presheaves over $X$ with values in $\emph{\textbf{K}}$, $u:\mathscr{F}_{1}\rightarrow\mathscr{F}_{1}$ a morphism (of presheaves), for each open subspace $U$ of $Y$ set $(\psi_{\ast}(u))_{U}=u_{\psi^{-1}(U)}$, then $\psi_{\ast}(u):\psi_{\ast}(\mathscr{F}_{1})
\rightarrow\psi_{\ast}(\mathscr{F}_{2})$ is a morphism (of presheaves), indeed $\psi_{\ast}:\psh_{\emph{\textbf{K}}}(X)
\rightarrow\psh_{\emph{\textbf{K}}}(Y)$, (also $\psi_{\ast}:\sh_{\emph{\textbf{K}}}(X)
\rightarrow\sh_{\emph{\textbf{K}}}(Y)$) is a covariant functor  where $\psh_{\emph{\textbf{K}}}(X)$ (resp. $\sh_{\emph{\textbf{K}}}(X)$) denotes the category of presheaves (resp. sheaves) over $X$ with values in the category
$\emph{\textbf{K}}$.\\
\end{bk}

\begin{bk}\label{I: 3.4.3} Let $Z$ be a third topological space, $\psi':Y\rightarrow Z$ a continuous map, and let $\psi''=\psi'\circ\psi$, then obviously we have $\psi''_{\ast}(\mathscr{F})=\psi'_{\ast}(\psi_{\ast}(\mathscr{F}))$ for each presheaf (resp. sheaf) $\mathscr{F}$ over $X$ with values in $\emph{\textbf{K}}$; moreover, for each morphism $u:\mathscr{F}\rightarrow\mathscr{G}$ of such presheaves, one has $\psi''_{\ast}(u)=\psi'_{\ast}(\psi_{\ast}(u))$. In other words, the functor $\psi''_{\ast}:\psh_{\emph{\textbf{K}}}(X)
\rightarrow\psh_{\emph{\textbf{K}}}(Z)$ is the composition of the functors $\psi_{\ast}:\psh_{\emph{\textbf{K}}}(X)
\rightarrow\psh_{\emph{\textbf{K}}}(Y)$ and $\psi'_{\ast}:\psh_{\emph{\textbf{K}}}(Y)
\rightarrow\psh_{\emph{\textbf{K}}}(Z)$, i.e. $\psi''_{\ast}=(\psi'\circ\psi)_{\ast}=\psi'_{\ast}\circ\psi_{\ast}$. \\
Moreover, for each open subspace $U$ of $Y$, the direct image of the presheaf $\mathscr{F}|_{\psi^{-1}(U)}$ by the restriction map $\psi|_{\psi^{-1}U}:\psi^{-1}(U)\rightarrow U$ is none other than the induced presheaf $(\psi_{\ast}(\mathscr{F}))|_{U}$, i.e. $(\psi|_{\psi^{-1}U})_{\ast}(\mathscr{F}|_{\psi^{-1}(U)})=
(\psi_{\ast}(\mathscr{F}))|_{U}$.  \\
\end{bk}

\begin{bk}\label{I: 3.4.4} Suppose that the category $\emph{\textbf{K}}$ admits the inductive limits, and let $\mathscr{F}$ be a presheaf over $X$ with values in $\emph{\textbf{K}}$, $\psi:X\rightarrow Y$ a continuous map and $x\in X$ a fixed point, then the canonical morphisms $\theta_{U}=\can_{\psi^{-1}(U),\mathscr{F}}:(\psi_{\ast}\mathscr{F})(U)=
\mathscr{F}(\psi^{-1}(U))\rightarrow\mathscr{F}_{x}$ where $U$ runs through the set of open neighborhoods of $\psi(x)$ in $Y$, form an inductive system, which gives, by passing the inductive limits, a morphism $\psi_{x,\mathscr{F}} :(\psi_{\ast}\mathscr{F})_{\psi(x)}\rightarrow\mathscr{F}_{x}$ of the fibers (this morphism, if there is no confusion, is also denoted by $\psi_{x}$ ), the morphism $\psi_{x,\mathscr{F}}$, in general, is neither a monomorphism nor an epimorphism.\\
In fact, $\psi_{x,-}:\Fib_{\psi(x)}\circ\psi_{\ast}\rightarrow\Fib_{x}$ is a natural transformation between the functors $$\Fib_{\psi(x)}\circ\psi_{\ast}:
\psh_{\emph{\textbf{K}}}(X)\rightarrow\emph{\textbf{K}}$$ and  $$\Fib_{x}:\psh_{\emph{\textbf{K}}}(X)\rightarrow\emph{\textbf{K}}$$ where $\Fib_{x}:\psh_{\emph{\textbf{K}}}(X)\rightarrow\emph{\textbf{K}}$ is the fiber functor at $x$ and similarly $\Fib_{\psi(x)}:\psh_{\emph{\textbf{K}}}(Y)\rightarrow\emph{\textbf{K}}$ is the fiber functor at $\psi(x)$. Because, if $u:\mathscr{F}_{1}\rightarrow\mathscr{F}_{2}$ is a morphism of presheaves over $X$ with values in $\emph{\textbf{K}}$, then for any open neighborhood $U$ of $\psi(x)$ in $Y$, one has $u_{x}\circ\psi_{x,\mathscr{F}_{1}}\circ\can'_{U,\mathscr{F}_{1}}=
u_{x}\circ\can_{\psi^{-1}(U),\mathscr{F}_{1}}=
\can_{\psi^{-1}(U),\mathscr{F}_{2}}\circ u_{\psi^{-1}(U)}$ where $\can_{\psi^{-1}(U),\mathscr{F}_{1}}:\mathscr{F}_{1}(\psi^{-1}(U))
\rightarrow(\mathscr{F}_{1})_{x}$ and $\can'_{U,\mathscr{F}_{1}}:(\psi_{\ast}\mathscr{F}_{1})(U)
\rightarrow(\psi_{\ast}\mathscr{F}_{1})_{\psi(x)}$ are the canonical morphisms (similarly for $\mathscr{F}_{2}$), moreover $\psi_{x,\mathscr{F}_{2}}\circ(\psi_{\ast}(u))_{\psi(x)}
\circ\can'_{U,\mathscr{F}_{1}}=\psi_{x,\mathscr{F}_{2}}
\circ\can'_{U,\mathscr{F}_{2}}\circ(\psi_{\ast}(u))_{U}=
\can_{\psi^{-1}(U),\mathscr{F}_{2}}\circ u_{\psi^{-1}(U)}$, hence by the universal property of the inductive limit $u_{x}\circ\psi_{x,\mathscr{F}_{1}}=
\psi_{x,\mathscr{F}_{2}}\circ(\psi_{\ast}(u))_{\psi(x)}$, therefore the following diagram is commutative $$\xymatrix{
(\psi_{\ast}\mathscr{F}_{1})_{\psi(x)}\ar[r]^{\psi_{x,\mathscr{F}_{1}}} \ar[d]^{(\psi_{\ast}(u))_{\psi(x)}} & (\mathscr{F}_{1})_{x} \ar[d]^{u_{x}} \\ (\psi_{\ast}\mathscr{F}_{2})_{\psi(x)}\ar[r]^{\psi_{x,\mathscr{F}_{2}} } &(\mathscr{F}_{2})_{x} .}$$
If $Z$ is a third topological space, $\psi':Y\rightarrow Z$ a continuous map, take $\psi''=\psi'\circ\psi$, then for any presheaf $\mathscr{F}$ over $X$ with values in $\emph{\textbf{K}}$, we have $\psi''_{x,\mathscr{F}}=
\psi_{x,\mathscr{F}}\circ\psi'_{\psi(x),\psi_{\ast}\mathscr{F}}:
(\psi''_{\ast}\mathscr{F})_{\psi''(x)}\rightarrow\mathscr{F}_{x}$. In other words, the natural transformations $\psi''_{x,-}:\Fib_{\psi''(x)}\circ\psi''_{\ast}\rightarrow\Fib_{x}$ and $\psi'_{\psi(x),-}\circ\psi_{\ast}:
\Fib_{\psi''(x)}\circ\psi''_{\ast}\rightarrow\Fib_{x}$ are equal (recall that if $\mu:F\rightarrow G$ is a natural transformation between the covariant functors $F,G:\mathscr{C}\rightarrow\mathscr{D}$, and let $H:\mathscr{A}\rightarrow\mathscr{C}$ be a third covariant functor, then the transformation $\mu\circ H:F\circ H\rightarrow G\circ H$ given by $(\mu\circ H)_{A}=\mu_{H(A)}$ for each object $A$ of $\mathscr{A}$, is a natural transformation). Because, for any open neighborhood $U$ of $\psi''(x)$ in $Z$, one has $\psi_{x,\mathscr{F}}\circ\psi'_{\psi(x),\psi_{\ast}\mathscr{F}}
\circ\can_{U,\psi''_{\ast}\mathscr{F}}=
\psi_{x,\mathscr{F}}\circ\can_{(\psi')^{-1}(U),\psi_{\ast}\mathscr{F}}=
\can_{(\psi'')^{-1}(U),\mathscr{F}}$ where
 $\can_{(\psi'')^{-1}(U),\mathscr{F}}:\mathscr{F}((\psi'')^{-1}(U))
 \rightarrow\mathscr{F}_{x}$, $\can_{(\psi')^{-1}(U),\psi_{\ast}\mathscr{F}}:
 (\psi_{\ast}\mathscr{F})((\psi')^{-1}(U))
 \rightarrow(\psi_{\ast}\mathscr{F})_{\psi(x)}$ and $\can_{U,\psi''_{\ast}\mathscr{F}}:(\psi''_{\ast}\mathscr{F})(U)
 \rightarrow(\psi''_{\ast}\mathscr{F})_{\psi''(x)}$ are the canonical morphisms, moreover $\psi''_{x,\mathscr{F}}\circ\can_{U,\psi''_{\ast}\mathscr{F}}=
 \can_{(\psi'')^{-1}(U),\mathscr{F}}$; hence by the universal property of the inductive limit $\psi''_{x,\mathscr{F}}=
\psi_{x,\mathscr{F}}\circ\psi'_{\psi(x),\psi_{\ast}\mathscr{F}}$.\\
\end{bk}

\begin{bk}\label{I: 3.4.5} Under the hypotheses \ref{I: 3.4.4}, suppose also that $\psi$ is an homeomorphism from $X$ onto its image, then
$\psi_{x,\mathscr{F}}:(\psi_{\ast}\mathscr{F})_{\psi(x)}
\rightarrow\mathscr{F}_{x}$ is an isomorphism. Because, since $\psi$ is injective and an open map, hence for any open neighborhood $V$ of $x$ in $X$, there exists (at least) an open subspace $U$ of $Y$ so that $V=\psi^{-1}(U)$ (or equivalently $\psi(V)=\psi(X)\cap U$), then take $\zeta_{V}=\can_{U,\psi_{\ast}\mathscr{F}}:\mathscr{F}(V)
\rightarrow(\psi_{\ast}\mathscr{F})_{\psi(x)}$ where  $\can_{U,\psi_{\ast}\mathscr{F}}
:(\psi_{\ast}\mathscr{F})(U)\rightarrow(\psi_{\ast}\mathscr{F})_{\psi(x)}$ is the canonical morphism. The morphisms $\zeta_{V}$ where $V$ runs through the set of open neighborhoods of $x$ in $X$, form an inductive system; because if $V,V'$ are such open neighborhoods in which $V'\subseteq V$, then we have $\zeta_{V'}\circ\rho_{V'}^{V}=
\can_{U',\psi_{\ast}\mathscr{F}}\circ\rho_{V'}^{V}=
\can_{U\cap U',\psi_{\ast}\mathscr{F}}\circ\omega_{U\cap U'}^{U'}\circ\rho_{V'}^{V}=\can_{U\cap U',\psi_{\ast}\mathscr{F}}\circ\rho_{\psi^{-1}(U\cap U')}^{\psi^{-1}( U')} \circ\rho_{V'}^{V}=\can_{U\cap U',\psi_{\ast}\mathscr{F}}\circ\rho_{V'}^{V}=\can_{U\cap U',\psi_{\ast}\mathscr{F}}\circ\omega_{\psi^{-1}(U\cap U')}^{\psi^{-1}( U)}=\can_{U,\psi_{\ast}\mathscr{F}}=\zeta_{V}$, hence by the universal property of the inductive limit there exists a (unique) morphism $\zeta_{x,\mathscr{F}}:
\mathscr{F}_{x}\rightarrow(\psi_{\ast}\mathscr{F})_{\psi(x)}$ so that $\zeta_{V}=\zeta_{x,\mathscr{F}}\circ\can_{V,\mathscr{F}}$ where $\can_{V,\mathscr{F}}:\mathscr{F}(V)\rightarrow\mathscr{F}_{x}$ is the canonical morphism. But, the definition of $\zeta_{V}$ is independent of choosing such $U$, because if there exists a second open subspace $W$ of $Y$ so that $V=\psi^{-1}(W)$, then $\can_{U,\psi_{\ast}\mathscr{F}}=\can_{U\cap W,\psi_{\ast}\mathscr{F}}\circ\omega_{U\cap W}^{U}=\can_{U\cap W,\psi_{\ast}\mathscr{F}}\circ\rho_{\psi^{-1}(U\cap W)}^{\psi^{-1}(U)}=\can_{U\cap W,\psi_{\ast}\mathscr{F}}=\can_{U\cap W,\psi_{\ast}\mathscr{F}}\circ\rho_{\psi^{-1}(U\cap W)}^{\psi^{-1}(W)}=\can_{U\cap W,\psi_{\ast}\mathscr{F}}\circ\omega_{U\cap W}^{W}=\can_{W,\psi_{\ast}\mathscr{F}}$, therefore by the universal property of the inductive limit,  $\zeta_{x,\mathscr{F}}\circ\psi_{x,\mathscr{F}}$ is the identity morphism. The morphism $\psi_{x,\mathscr{F}}\circ\zeta_{x,\mathscr{F}}$, using the universal property of the inductive limit once more, is also the identity morphism. \\ This applies in particular to the injection map $j:X\rightarrow Y$ where $X$ is a subspace of $Y$. \\
\end{bk}

\begin{bk}\label{I: 3.4.6} Suppose $\emph{\textbf{K}}$ is the category of groups or the category of rings. Note that in this case, $\emph{\textbf{K}}$ admits filtered (directed) inductive limits.  If $\mathscr{F}$ is a sheaf over $X$ with values in $\emph{\textbf{K}}$, of support $S$, and $\psi:X\rightarrow Y$ a continuous map, then $\Supp(\psi_{\ast}\mathscr{F})\subseteq\overline{\psi(S)}$ where $\overline{\psi(S)}$ is the closure of $\psi(S)$ in $Y$. Because, suppose on the contrary and then take $y\in\Supp(\psi_{\ast}\mathscr{F})\setminus\overline{\psi(S)}$, then there exists an open neighborhood $V$ of $y$ in $Y$ so that $V\cap\psi(S)=\emptyset$, and also there exists an (nonzero) element $s\in(\psi_{\ast}\mathscr{F})(V)$ in which the germ of $s$ at $y$ is nonzero, i.e. $0\neq s_{y}\in(\psi_{\ast}\mathscr{F})_{y}$. Since $\mathscr{F}$ is a sheaf therefore $\psi^{-1}(V)\neq\emptyset$. On the other hand, since $s\in(\psi_{\ast}\mathscr{F})(V)=\mathscr{F}(\psi^{-1}(V))$ is a nonzero element therefore there exists some point $x\in\psi^{-1}(V)$ in which $s_{x}\in\mathscr{F}_{x}$ is nonzero (because $\mathscr{F}$ is a sheaf in the sense of \ref{1.2}), hence $x\in\Supp(\mathscr{F})=S$ a contradiction.
Under the same hypotheses, if $X$ is a subspace of $Y$ and $j:X\rightarrow Y$ is the injection map, then the sheaf $j_{\ast}\mathscr{F}$ induces the sheaf $\mathscr{F}$ over $X$ (i.e. for each $x\in X$, $(j_{\ast}\mathscr{F})_{x}=\mathscr{F}_{x}$), if moreover $X$ is a closed subspace of $Y$, then the sheaf $j_{\ast}\mathscr{F}$ induces the sheaf $\mathscr{F}$ over $X$ and the zero sheaf over $Y-X$ (i.e. for each $y\in Y$, if $y\in X$ then $(j_{\ast}\mathscr{F})_{y}=\mathscr{F}_{y}$ and if $y\in Y-X$, then $(j_{\ast}\mathscr{F})_{y}=\{0\}$).\\
\end{bk}

\section{Inverse image of presheaves}

\begin{bk}\label{I: 3.5.1} Under the hypotheses of \ref{I: 3.4.1}, if $\mathscr{F}$ (resp. $\mathscr{G}$) is a presheaf over $X$ (resp. $Y$) with values in a category $\emph{\textbf{K}}$, each morphism $u:\mathscr{G}\rightarrow\psi_{\ast}(\mathscr{F})$ of presheaves over $Y$ is also called a $\psi-$morphism from $\mathscr{G}$ into $\mathscr{F}$, and by abuse of notation it is denoted by $\mathscr{G}\rightarrow\mathscr{F}$. The set of $\psi-$morphisms from $\mathscr{G}$ into $\mathscr{F}$, $\Hom_{Y}(\mathscr{G},\psi_{\ast}(\mathscr{F}))$, is also denoted by $\Hom_{\psi}(\mathscr{G},\mathscr{F})$. For each couple $(U,V)$, where $U$ is an open subspace of $X$, $V$ an open subspace of $Y$ with $\psi(U)\subseteq V$ (or equivalently $U\subseteq\psi^{-1}(V)$), one has a morphism $u_{U,V}:\mathscr{G}(V)\rightarrow\mathscr{F}(U)$ which is, in fact, the composition of the morphism $u_{V}:\mathscr{G}(V)\rightarrow(\psi_{\ast}\mathscr{F})(V)=
\mathscr{F}(\psi^{-1}(V))$ with the restriction morphism $\rho_{U}^{\psi^{-1}(V)}:\mathscr{F}(\psi^{-1}(V))\rightarrow
\mathscr{F}(U)$, i.e. $u_{U,V}=\rho_{U}^{\psi^{-1}(V)}\circ u_{V}$. Obviously, for each such two couples $(U,V)$ and $(U',V')$ with $U'\subseteq U$ and $V'\subseteq V$, the following diagram is commutative $$\xymatrix{
\mathscr{G}(V) \ar[r]^{u_{U,V}} \ar[d]^{\varrho_{V'}^{V}} & \mathscr{F}(U) \ar[d]^{\rho_{U'}^{U}} \\ \mathscr{G}(V')\ar[r]^{u_{U',V'}} & \mathscr{F}(U') .}$$ Conversely, given a family of morphisms $(u_{U,V}:\mathscr{G}(V)\rightarrow\mathscr{F}(U))$ where $V$ is an arbitrary open subspace of $Y$, $U$ an arbitrary open subspace of $X$ provided that $\psi(U)\subseteq V$, and for each such two couples $(U,V)$ and $(U',V')$ with $U'\subseteq U$ and $V'\subseteq V$, the above diagram is commutative. Then this family defines a $\psi-$morphism $u:\mathscr{G}\rightarrow\mathscr{F}$ by taking $u_{V}=u_{\psi^{-1}(V),V}$ for each open subspace $V$ of $Y$.\\
Suppose that the category $\emph{\textbf{K}}$ admits the projective limits, $\mathscr{F}$ (resp. $\mathscr{G}$ ) is a sheaf over $X$ (resp. $Y$) with values in  $\emph{\textbf{K}}$, $\mathfrak{B}$ and $\mathfrak{B}'$ are bases for the topologies of $X$ and $Y$ respectively, also suppose that a family $(u_{U,V}:\mathscr{G}(V)\rightarrow\mathscr{F}(U))$ of morphisms is given where $V$ is an arbitrary open subspace belongs to $\mathfrak{B}'$, $U$ an arbitrary open subspace belongs to $\mathfrak{B}$ provided that $\psi(U)\subseteq V$, and for each such two couples $(U,V)$ and $(U',V')$ with $U'\subseteq U$ and $V'\subseteq V$, the above diagram is commutative. Then this family also defines a $\psi-$morphism $u:\mathscr{G}\rightarrow\mathscr{F}$. Because, for each fixed open subspace $W$ of $Y$, first we define a morphism $u_{W}:\mathscr{G}(W)\rightarrow(\psi_{\ast}\mathscr{F})(W)$ as follows.
For each fixed $V\in\mathfrak{B}'$ with $V\subseteq W$, the family of morphisms $(u_{U,V})_{U}$, where $U$ runs through the set of elements of $\mathfrak{B}$ contained in $\psi^{-1}(V)$, constitute a projective system and so by the universal property of the projective limit there exists a unique morphism $\theta_{V}:\mathscr{G}(V)\rightarrow(\psi_{\ast}\mathscr{F})(V)=
\mathscr{F}(\psi^{-1}(V))=
\lim\limits_{U\in\mathfrak{B}}\mathscr{F}(U)$ in which for each
$U\in\mathfrak{B}$ with $U\subseteq\psi^{-1}(V)$, one has $u_{U,V}=\can_{U,\psi^{-1}(V)}\circ\theta_{V}$ where $\can_{U,\psi^{-1}(V)}:\mathscr{F}(\psi^{-1}(V))\rightarrow\mathscr{F}(U)$ is the canonical morphism, then take $\zeta_{V}=\theta_{V}\circ\varrho_{V}^{W}:
\mathscr{G}(W)\rightarrow(\psi_{\ast}\mathscr{F})(V)$. The family of morphisms $(\zeta_{V})$, where $V$ runs through the set of elements of $\mathfrak{B}'$ contained in $W$, constitute a projective system and by the universal property of the projective limit there exists a unique morphism $u_{W}:\mathscr{G}(W)\rightarrow(\psi_{\ast}\mathscr{F})(W)=
\lim\limits_{V\in\mathfrak{B}'}(\psi_{\ast}\mathscr{F})(V)$ in which for each $V\in\mathfrak{B}'$ with $V\subseteq W$,
$\zeta_{V}=\can'_{V,W}\circ u_{W}$ where $\can'_{V,W}:(\psi_{\ast}\mathscr{F})(W)
\rightarrow(\psi_{\ast}\mathscr{F})(V)$ is the canonical morphism. Then it is easy to observe that
 $u:\mathscr{G}\rightarrow\psi_{\ast}\mathscr{F}$ is actually a morphism of
sheaves over $Y$ (or equivalently $u$ is a $\psi-$morphism).\\
Suppose that the category $\emph{\textbf{K}}$ admits the inductive limits, $u:\mathscr{G}\rightarrow\mathscr{F}$ a $\psi-$morphism and $x\in X$ a fixed point. The family of morphisms $(\can_{\psi^{-1}(V),\mathscr{F}}\circ u_{V}:\mathscr{G}(V)\rightarrow\mathscr{F}_{x})$ where $\can_{\psi^{-1}(V),\mathscr{F}}:\mathscr{F}(\psi^{-1}(V))
\rightarrow\mathscr{F}_{x}$ is the canonical morphism and $V$ runs through the set of open neighborhoods of $\psi(x)$ in $Y$, constitute an inductive system which gives, by passing to the limits, a unique morphism $u_{x}:\mathscr{G}_{\psi(x)}\rightarrow\mathscr{F}_{x}$ so that for each such $V$, $\can_{\psi^{-1}(V),\mathscr{F}}\circ u_{V}=u_{x}\circ\can_{V,\mathscr{G}}$ where $\can_{V,\mathscr{G}}:\mathscr{G}(V)\rightarrow\mathscr{G}_{\psi(x)}$ is the canonical morphism. In fact, it is easy to observe that  $u_{x}=\psi_{x,\mathscr{F}}\circ u_{\psi(x)}$, \ref{I: 3.4.4}.\\
\end{bk}

\begin{bk}\label{I: 3.5.2} Let $\emph{\textbf{K}}$ be a category, then we can form a category over $\emph{\textbf{K}}$ whose objects are the couples $(X,\mathscr{F})$ where $X$ is a topological space and $\mathscr{F}$ is a presheaf over $X$ with values in $\emph{\textbf{K}}$, and the morphisms are the pairs $(\psi, u):(X,\mathscr{F})\rightarrow(Y, \mathscr{G})$ where $\psi:X\rightarrow Y$ is a continuous map and $u:\mathscr{G}\rightarrow\mathscr{F}$ a $\psi-$morphism.\\
\end{bk}

\begin{bk}\label{I: 3.5.3} Let $\psi:X\rightarrow Y$ be a continuous map, $\mathscr{G}$ a presheaf over $Y$ with values in a category $\emph{\textbf{K}}$. The pair $(\mathscr{G}', \rho)$ where $\mathscr{G}'$ is a sheaf over $X$ with values in $\emph{\textbf{K}}$ and $\rho:\mathscr{G}\rightarrow\mathscr{G}'$ a $\psi-$morphism (in other words a morphism $\mathscr{G}\rightarrow\psi_{\ast}(\mathscr{G}')$ of presheaves over $Y$) is said to be the \emph{inverse} \emph{image} \emph{of} $\mathscr{G}$  \emph{by} $\psi$, if for each sheaf $\mathscr{F}$ over $X$ with values in  $\emph{\textbf{K}}$, the function $$\Hom_{X}(\mathscr{G}',\mathscr{F})
\rightarrow\Hom_{\psi}(\mathscr{G},\mathscr{F})=
\Hom_{Y}(\mathscr{G},\psi_{\ast}(\mathscr{F}))$$ which transforms each $\nu$ into $\psi_{\ast}(\nu)\circ\rho$ is a bijection. This function is functorial at $\mathscr{F}$; more precisely, for each morphism of sheaves $u:\mathscr{F}_{1}\rightarrow\mathscr{F}_{2}$ over $X$, the following diagram is commutative $$\xymatrix{
\Hom_{X}(\mathscr{G}',\mathscr{F}_{1})\ar[r]^{} \ar[d]^{} & \Hom_{Y}(\mathscr{G},\psi_{\ast}\mathscr{F}_{1}) \ar[d]^{} \\ \Hom_{X}(\mathscr{G}',\mathscr{F}_{2})\ar[r]^{} & \Hom_{Y}(\mathscr{G},\psi_{\ast}\mathscr{F}_{2}) .}$$ Since the pair $(\mathscr{G}', \rho)$ is a solution of the universal problem, therefore, if it exists, it is unique up to a unique isomorphism. We will denote $\mathscr{G}'$ by $\psi^{\ast}(\mathscr{G})$ and $\rho$ by $\rho_{\mathcal{G}}$, i.e. $\psi^{\ast}(\mathscr{G})=\mathscr{G}'$ and $\rho_{\mathcal{G}}=\rho$. By abuse of the notation, we shall  say that $\psi^{\ast}(\mathscr{G})$ is the inverse image sheaf of $\mathscr{G}$ by $\psi$ (i.e. we shall purposely drop the morphism  $\rho_{\mathcal{G}}$ and do not mention it when we are talking about the inverse image sheaf  $\psi^{\ast}(\mathscr{G})$) provided that $\psi^{\ast}(\mathscr{G})$ is equipped with a canonical $\psi-$morphism $\rho_{\mathcal{G}}:\mathscr{G}\rightarrow\psi^{\ast}(\mathscr{G})$, or more precisely the canonical homomorphism of presheaves over $Y$:
$$\rho_{\mathcal{G}}:\mathscr{G}\rightarrow
\psi_{\ast}(\psi^{\ast}(\mathscr{G})) .$$ For each morphism of sheaves $\nu:\psi^{\ast}(\mathscr{G})\rightarrow\mathscr{F}$ (where $\mathscr{F}$ is a sheaf over $X$ with values in  $\emph{\textbf{K}}$), we set $\nu^{\flat}=\psi_{\ast}(\nu)\circ\rho_{\mathcal{G}}$.  By definition, each morphism $u:\mathscr{G}\rightarrow\psi_{\ast}(\mathscr{F})$ is of the form $\nu^{\flat}$ for one $\nu$ and only one, which we shall denote it by $u^{\sharp}$, i.e. $u^{\sharp}=\nu$. In other words, each morphism of presheaves $u:\mathscr{G}\rightarrow\psi_{\ast}(\mathscr{F})$ uniquely decomposes as $\xymatrix{u:\mathscr{G}\ar[r]^{\rho_{\mathcal{G}}} & \psi_{\ast}(\psi^{\ast}(\mathscr{G})) \ar[r]^{\psi_{\ast}(u^{\sharp})} & \psi_{\ast}(\mathscr{F}) .}$ \\
\end{bk}

\begin{bk}\label{I: 3.5.4} Now assume that the category $\emph{\textbf{K}}$ is such that each presheaf $\mathscr{G}$ over $Y$ with values in $\emph{\textbf{K}}$ admits an inverse image sheaf by $\psi$, which we denote it by $\psi^{\ast}(\mathscr{G})$. We shall observe that $\psi^{\ast}$ is actually a (covariant) functor from the category $\psh_{\emph{\textbf{K}}}(Y)$ into the category  $\sh_{\emph{\textbf{K}}}(X)$ in which the bijection $\nu\rightsquigarrow\nu^{\flat}$ is an isomorphism between the following bifunctors $$\Hom_{X}(\psi^{\ast}(-),-):
\psh_{\emph{\textbf{K}}}(Y)\times\sh_{\emph{\textbf{K}}}(X)
\rightarrow\sets$$ and $$\Hom_{Y}(-,\psi_{\ast}(-)):
\psh_{\emph{\textbf{K}}}(Y)\times\sh_{\emph{\textbf{K}}}(X)
\rightarrow\sets .$$ Because, for an arbitrary morphism $u:\mathscr{G}\rightarrow\mathscr{H}$ of presheaves over $Y$, we define $\psi^{\ast}(u)=(\rho_{\mathcal{\mathcal{H}}}\circ u)^{\sharp}:\psi^{\ast}(\mathscr{G})\rightarrow\psi^{\ast}(\mathscr{H})$ where $\rho_{\mathcal{H}}:\mathscr{H}\rightarrow
\psi_{\ast}(\psi^{\ast}(\mathscr{H}))$ is the canonical morphism. If $\Id_{\mathscr{G}}:\mathscr{G}\rightarrow\mathscr{G}$ is the identity morphism, then by the definition $\psi^{\ast}(\Id_{\mathscr{G}})=
(\rho_{\mathcal{G}}\circ\Id_{\mathscr{G}})^{\sharp}=
(\rho_{\mathcal{G}})^{\sharp}=\Id_{\psi^{\ast}(\mathscr{G})}$ (since the map $\Hom_{X}(\psi^{\ast}(\mathscr{G}),\psi^{\ast}(\mathscr{G}))
\rightarrow\Hom_{Y}(\mathscr{G},\psi_{\ast}(\psi^{\ast}(\mathscr{G})))$ given by $\nu\rightsquigarrow\nu^{\flat}=\psi_{\ast}(\nu)\circ\rho_{\mathcal{G}}$ is bijective and $\psi_{\ast}$ is a functor therefore $(\rho_{\mathcal{G}})^{\sharp}=\Id_{\psi^{\ast}(\mathscr{G})}$ ). Now let $u:\mathscr{G}\rightarrow\mathscr{H}$ and $v:\mathscr{H}\rightarrow\mathscr{L}$ be arbitrary morphisms of presheaves over $Y$, we have then $\psi^{\ast}(v\circ u)=\psi^{\ast}(v)\circ\psi^{\ast}(u)$,
because $\psi_{\ast}((\rho_{\mathcal{L}}\circ v)^{\sharp}\circ(\rho_{\mathcal{H}}\circ u)^{\sharp})\circ\rho_{\mathcal{G}}=
\psi_{\ast}((\rho_{\mathcal{L}}\circ v)^{\sharp})\circ\psi_{\ast}((\rho_{\mathcal{H}}\circ u)^{\sharp})\circ\rho_{\mathcal{G}}=\psi_{\ast}((\rho_{\mathcal{L}}\circ v)^{\sharp})\circ\rho_{\mathcal{H}}\circ u=\rho_{\mathcal{L}}\circ v\circ u=\psi_{\ast}((\rho_{\mathcal{L}}\circ v\circ u)^{\sharp})\circ\rho_{\mathcal{G}}$, therefore $(\rho_{\mathcal{L}}\circ v\circ u)^{\sharp}=(\rho_{\mathcal{L}}\circ v)^{\sharp}\circ(\rho_{\mathcal{H}}\circ u)^{\sharp}$, in other words $\psi^{\ast}(v\circ u)=\psi^{\ast}(v)\circ\psi^{\ast}(u)$. Hence, $\psi^{\ast}:\psh_{\emph{\textbf{K}}}(Y)
\rightarrow\sh_{\emph{\textbf{K}}}(X)$ is actually a covariant functor.\\
For each sheaf $\mathscr{F}$ over $X$ with values in $\emph{\textbf{K}}$, let $i_{\mathscr{F}}$ to be the identity morphism of $\psi_{\ast}(\mathscr{F})$, i.e. $i_{\mathscr{F}}=\Id_{\psi_{\ast}(\mathscr{F})}:
\psi_{\ast}(\mathscr{F})
\rightarrow\psi_{\ast}(\mathscr{F})$ and take $\sigma_{\mathscr{F}}=(i_{\mathscr{F}})^{\sharp}:
\psi^{\ast}(\psi_{\ast}(\mathscr{F}))\rightarrow\mathscr{F}$. Then for each morphism $u:\mathscr{G}\rightarrow\psi_{\ast}(\mathscr{F})$ of presheaves over $Y$, we have the decomposition $$\xymatrix{u^{\sharp}:\psi^{\ast}(\mathscr{G})\ar[r]^{\psi^{\ast}(u)} & \psi^{\ast}(\psi_{\ast}(\mathscr{F})) \ar[r]^{\sigma_{\mathscr{F}}} & \mathscr{F} .}$$ Because, one has $\psi_{\ast}(\sigma_{\mathscr{F}}\circ\psi^{\ast}(u))
\circ\rho_{\mathcal{G}}=
\psi_{\ast}(\sigma_{\mathscr{F}})\circ\psi_{\ast}(\psi^{\ast}(u))
\circ\rho_{\mathcal{G}}=\psi_{\ast}((i_{\mathscr{F}})^{\sharp})
\circ\psi_{\ast}((\rho_{\psi_{\ast}(\mathscr{F})}\circ u)^{\sharp})\circ\rho_{\mathcal{G}}=
\psi_{\ast}((\Id_{\psi_{\ast}(\mathscr{F})})^{\sharp})
\circ\rho_{\psi_{\ast}(\mathscr{F})}\circ u=\Id_{\psi_{\ast}(\mathscr{F})}\circ u=u$, therefore $u^{\sharp}=\sigma_{\mathscr{F}}\circ\psi^{\ast}(u)$. We say that the morphism $\sigma_{\mathscr{F}}$ is canonical. \\
\end{bk}

\begin{bk}\label{I: 3.5.5} Let $\psi':Y\rightarrow Z$ be a continuous map, suppose that each presheaf $\mathscr{H}$ over $Z$ with values in $\emph{\textbf{K}}$ admits an inverse image $\psi'^{\ast}(\mathscr{H})$ by $\psi'$. Then (under the hypotheses of \ref{I: 3.5.4}) each presheaf $\mathscr{H}$ over $Z$ with values in $\emph{\textbf{K}}$ admits an inverse  image by $\psi''=\psi'\circ\psi$ and we have a canonical isomorphism $\psi''^{\ast}\simeq\psi^{\ast}\circ\psi'^{\ast}$
between the functors $$\psi''^{\ast}, \psi^{\ast}\circ\psi'^{\ast}:\psh_{\emph{\textbf{K}}}(Z)\rightarrow
\sh_{\emph{\textbf{K}}}(X) .$$ First we show that the pair $(\mathscr{H}',\rho)$ including the sheaf $\mathscr{H}'=\psi^{\ast}(\psi'^{\ast}(\mathscr{H}))$ together with the morphism $$\rho=\psi'_{\ast}(\rho_{\psi'^{\ast}(\mathscr{H})})\circ\rho_{\mathcal{H}}:
\mathscr{H}\rightarrow\psi''_{\ast}(\mathscr{H}')$$ is an inverse image of the presheaf $\mathscr{H}$ by the map $\psi'':X\rightarrow Z$, where $\rho_{\mathcal{H}}:\mathscr{H}
\rightarrow\psi'_{\ast}(\psi'^{\ast}(\mathscr{H}))$ and $\rho_{\psi'^{\ast}(\mathscr{H})}:\psi'^{\ast}(\mathscr{H})
\rightarrow\psi_{\ast}(\psi^{\ast}(\psi'^{\ast}(\mathscr{H})))$ are the canonical morphisms. Let $\mathscr{F}$ be a sheaf over $X$ with values in $\emph{\textbf{K}}$, it is sufficient to show that the map $\theta:\Hom_{X}(\mathscr{H}',\mathscr{F})\rightarrow
\Hom_{Z}(\mathscr{H},\psi''_{\ast}(\mathscr{F}))$ given by
$\nu\rightsquigarrow\psi''_{\ast}(\nu)\circ\rho$ is bijective. But we have $\theta=\theta_{2}\circ\theta_{1}$ where
the maps $$\theta_{1}:\Hom_{X}(\psi^{\ast}(\psi'^{\ast}(\mathscr{H})),
\mathscr{F})\rightarrow
\Hom_{Y}(\psi'^{\ast}(\mathscr{H}),\psi_{\ast}(\mathscr{F}))$$ $$\nu\rightsquigarrow\psi_{\ast}(\nu)\circ\rho_{\psi'^{\ast}(\mathscr{H})}$$
and $$\theta_{2}:\Hom_{Y}(\psi'^{\ast}(\mathscr{H}),\psi_{\ast}(\mathscr{F}))
\rightarrow\Hom_{Z}(\mathscr{H},\psi'_{\ast}(\psi_{\ast}(\mathscr{F})))$$
$$\nu\rightsquigarrow\psi'_{\ast}(\nu)\circ\rho_{\mathcal{H}}$$ are bijective. Instead of proving the second part of the assertion, we prove a more general assertion which in particular gives up our desired assertion. Suppose that any presheaf over $Y$ with values in $\emph{\textbf{K}}$ admits an inverse image sheaf by (the continuous map) $\psi:X\rightarrow Y$; suppose also that for each presheaf $\mathscr{G}$ over $Y$ with values in $\emph{\textbf{K}}$, the pairs $(\sigma^{\ast}(\mathscr{G}),\rho_{\sigma^{\ast}(\mathscr{G})})$ and $(\eta^{\ast}(\mathscr{G}),\rho_{\eta^{\ast}(\mathscr{G})})$ are inverse images of $\mathscr{G}$ by $\psi$. In fact, according to \ref{I: 3.5.4}, $\sigma^{\ast},\eta^{\ast}:\psh_{\emph{\textbf{K}}}(Y)
\rightarrow\sh_{\emph{\textbf{K}}}(X)$ are covariant functors. Also, by the universal property of the inverse image sheaf, for each presheaf $\mathscr{G}$ over $Y$ with values in $\emph{\textbf{K}}$, there exists a unique isomorphism (of sheaves) $\zeta_{\mathscr{G}}:\sigma^{\ast}(\mathscr{G})
\rightarrow\eta^{\ast}(\mathscr{G})$ so that $\rho_{\eta^{\ast}(\mathscr{G})}=\psi_{\ast}(\zeta_{\mathscr{G}})
\circ\rho_{\sigma^{\ast}(\mathscr{G})}$. We claim that $\zeta:\sigma^{\ast}\rightarrow\eta^{\ast}$ is a natural transformation (in fact an isomorphism) of functors. Let $\theta:\mathscr{G}\rightarrow\mathscr{H}$ be an arbitrary morphism of presheaves over $Y$ with values in $\emph{\textbf{K}}$. One has $\psi_{\ast}(\eta^{\ast}(\theta)\circ\zeta_{\mathscr{G}})
\circ\rho_{\sigma^{\ast}(\mathscr{G})}=
\psi_{\ast}((\rho_{\eta^{\ast}(\mathscr{H})}
\circ\theta)^{\sharp})\circ\psi_{\ast}(\zeta_{\mathscr{G}})
\circ\rho_{\sigma^{\ast}(\mathscr{G})}=
\psi_{\ast}((\rho_{\eta^{\ast}(\mathscr{H})}
\circ\theta)^{\sharp})\circ\rho_{\eta^{\ast}(\mathscr{G})}=
\rho_{\eta^{\ast}(\mathscr{H})}\circ\theta=
\psi_{\ast}(\zeta_{\mathscr{H}}
\circ\sigma^{\ast}(\theta))\circ\rho_{\sigma^{\ast}(\mathscr{G})}$, but the map $$\Hom_{X}(\sigma^{\ast}(\mathscr{G}),\eta^{\ast}(\mathscr{H}))
\rightarrow\Hom_{Y}(\mathscr{G},\psi_{\ast}(\eta^{\ast}(\mathscr{H})))$$ given by $\nu\rightsquigarrow
\psi_{\ast}(\nu)\circ\rho_{\sigma^{\ast}(\mathscr{G})}$ is bijective, therefore $\eta^{\ast}(\theta)\circ\zeta_{\mathscr{G}}=\zeta_{\mathscr{H}}
\circ\sigma^{\ast}(\theta)$. In other words, the following diagram is commutative $$\xymatrix{
\sigma^{\ast}(\mathscr{G}) \ar[r]^{\zeta_{\mathscr{G}}} \ar[d]^{\sigma^{\ast}(\theta)} & \eta^{\ast}(\mathscr{G}) \ar[d]^{\eta^{\ast}(\theta)} \\ \sigma^{\ast}(\mathscr{H})\ar[r]^{\zeta_{\mathscr{H}}} &\eta^{\ast}(\mathscr{H}) .}$$ \\
\end{bk}

\begin{bk}\label{I: 3.5.6} Take in particular for $\psi$ the identity map $1_{X}:X\rightarrow X$. If $1_{X}^{\ast}(\mathscr{G})$, the inverse image sheaf by $\psi$ of a presheaf $\mathscr{G}$ over $X$ with values in $\emph{\textbf{K}}$, exists, then we call $1_{X}^{\ast}(\mathscr{G})$ the sheaf associated to the presheaf $\mathscr{G}$. Obviously, in this case, each morphism $u:\mathscr{G}\rightarrow\psi_{\ast}(\mathscr{F})=
(1_{X})_{\ast}(\mathscr{F})=\mathscr{F}$ from $\mathscr{G}$ into the sheaf $\mathscr{F}$ ($\mathscr{F}$ is a sheaf over $X$ with values in $\emph{\textbf{K}}$) decomposes in the unique way as $\xymatrix{u:\mathscr{G} \ar[r]^{\rho_{\mathcal{G}}} & 1_{X}^{\ast}(\mathscr{G}) \ar[r]^{u^{\sharp}} & \mathscr{F} .}$ \\
\end{bk}

\section{Simple and locally simple sheaves}

\begin{bk}\label{I: 3.6.1} We say that a presheaf $\mathscr{F}$ over $X$ with values in $\emph{\textbf{K}}$, is constant if for each nonempty open subspace $U$ of $X$ the canonical morphism $\rho_{U}^{X}:\mathscr{F}(X)\rightarrow\mathscr{F}(U)$ is an isomorphism; note that $\mathscr{F}$ is not necessarily a sheaf. We say that a sheaf $\mathscr{F}$ (over $X$ with values in $\emph{\textbf{K}}$) is simple if it is  associated to a constant presheaf over $X$ with values in $\emph{\textbf{K}}$, \ref{I: 3.5.6}. We say that a sheaf $\mathscr{F}$ (over $X$ with values in $\emph{\textbf{K}}$) is locally simple if each point $x\in X$ admits an open neighborhood $U$ (in $X$) such that the sheaf $\mathscr{F}|_{U}$ is simple.\\
\end{bk}

\begin{bk}\label{I: 3.6.2} Suppose that $X$ is an irreducible topological space; then the following conditions are equivalent. \\
$\leta)$ $\mathscr{F}$ is a constant presheaf over $X$.\\
$\letb)$ $\mathscr{F}$ is a simple sheaf over $X$.\\
$\letc)$ $\mathscr{F}$ is a locally simple sheaf over $X$.\\

$\leta)\Rightarrow\letb)$: Let $\mathscr{F}$ be a constant presheaf over $X$. We shall prove that $\mathscr{F}$ is a sheaf (over $X$ with values in $\emph{\textbf{K}}$). Suppose that $U$ is an open subspace of $X$, let $(U_{\alpha})_{\alpha\in I}$ be a covering for $U$ consisted of open subspaces contained in $U$, $T$ an arbitrary object of $\emph{\textbf{K}}$, then the function $$\Hom_{\emph{\textbf{K}}}(T,\mathscr{F}(U))
\rightarrow\prod\limits_{\alpha\in I} \Hom_{\emph{\textbf{K}}}(T,\mathscr{F}(U_{\alpha}))$$ given by $f\rightsquigarrow(\rho_{\alpha}\circ f)$ is a bijective from $\Hom_{\emph{\textbf{K}}}(T,\mathscr{F}(U))$ over the set of all sequences $(f_{\alpha})\in\prod\limits_{\alpha\in I} \Hom_{\emph{\textbf{K}}}(T,\mathscr{F}(U_{\alpha}))$ so that $\rho_{\alpha\beta}\circ f_{\alpha}=\rho_{\beta\alpha}\circ f_{\beta}$ for each pair of indices $(\alpha,\beta)$.\\
Case 1: Suppose that $U\neq\emptyset$. Also suppose that for each $\alpha$, $\rho_{\alpha}\circ f=\rho_{\alpha}\circ g$ for some morphisms $f, g\in\Hom_{\emph{\textbf{K}}}(T,\mathscr{F}(U))$. Since $U$ is nonempty there exists some $\beta\in I$ so that $U_{\beta}\neq\emptyset$ and so by the hypothesis $\rho_{\beta}=\rho_{U_{\beta}}^{X}\circ(\rho_{U}^{X})^{-1}$ is an isomorphism (hence a monomorphism) therefore, from $\rho_{\beta}\circ f=\rho_{\beta}\circ g$, we conclude that $f=g$. Now let $(f_{\alpha})\in\prod\limits_{\alpha\in I} \Hom_{\emph{\textbf{K}}}(T,\mathscr{F}(U_{\alpha}))$ be a family of morphisms so that $\rho_{\alpha\beta}\circ f_{\alpha}=\rho_{\beta\alpha}\circ f_{\beta}$ for each pair of indices $(\alpha,\beta)$. Take $f=(\rho_{\beta})^{-1}\circ f_{\beta}$ where $U_{\beta}\neq\emptyset$. Now for each $\alpha\in I$, one can observe that $\rho_{\alpha}\circ f=f_{\alpha}$. Because, if $U_{\alpha}\neq\emptyset$ then $U_{\alpha}\cap U_{\beta}\neq\emptyset$ since $X$ is irreducible, and we have $(\rho_{\alpha})^{-1}\circ f_{\alpha}=(\rho_{U_{\alpha}\cap U_{\beta}}^{U})^{-1}\circ\rho_{\alpha\beta}\circ f_{\alpha}=(\rho_{U_{\alpha}\cap U_{\beta}}^{U})^{-1}\circ\rho_{\beta\alpha}\circ f_{\beta}=(\rho_{\beta})^{-1}\circ f_{\beta}$ (this, in particular, also tells us that the definition of $f$ is independent of choosing such  $U_{\beta}$). But, if $U_{\alpha}=\emptyset$, then $U_{\alpha}\cap U_{\beta}=\emptyset=U_{\alpha}$ and we have $f_{\alpha}=\rho_{\alpha\beta}\circ f_{\alpha}=\rho_{\beta\alpha}\circ f_{\beta}=\rho_{U_{\alpha}\cap U_{\beta}}^{U}\circ\rho_{\beta}^{-1}\circ f_{\beta}=\rho_{\alpha}\circ f$.\\ Case 2: $U=\emptyset$. Then for each $\alpha\in I$, $U_{\alpha}=\emptyset=U$, hence $\rho_{\alpha}$ and $\rho_{\alpha\beta}$ are the identity morphisms, therefore the assertion easily implies.\\
Since $\mathscr{F}$ is a sheaf hence $\mathscr{F}=1_{X}^{\ast}(\mathscr{F})$, therefore $\mathscr{F}$ is a simple sheaf.  \\
$ \letc)\Rightarrow\leta)$: Suppose that $\mathscr{F}$ is a locally simple sheaf over $X$, hence there exists a covering $(U_{\alpha})$ for $X$ consisted of nonempty open subspaces so that for each $\alpha$, $\mathscr{F}|_{U_{\alpha}}$ is a simple sheaf. Therefore for each $\alpha$, there exists a constant presheaf $\mathscr{G}_{\alpha}$ over $U_{\alpha}$ (with values in $\emph{\textbf{K}}$) so that $1^{\ast}_{U_{\alpha}}(\mathscr{G}_{\alpha})=\mathscr{F}|_{U_{\alpha}}$.
Since $U_{\alpha}$ is irreducible hence by the above argument, i.e. $\leta)\Rightarrow\letb)$, $\mathscr{G}_{\alpha}$ is a sheaf therefore $\mathscr{G}_{\alpha}=1^{\ast}_{U_{\alpha}}(\mathscr{G}_{\alpha})$, and so $\mathscr{F}|_{U_{\alpha}}$ is a constant sheaf. Now one can observe that $\mathscr{F}$ is indeed a constant sheaf; because take a nonempty open subspace $U$ of $X$ and then for each $\alpha$ set $f_{\alpha}=(\rho_{U\cap U_{\alpha}}^{U_{\alpha}})^{-1}\circ\rho_{U\cap U_{\alpha}}^{U}:\mathscr{F}(U)\rightarrow\mathscr{F}(U_{\alpha})$. The family of morphisms $(f_{\alpha})\in\prod\limits_{\alpha}
\Hom(\mathscr{F}(U),\mathscr{F}(U_{\alpha}))$ obviously satisfies in the condition $\rho_{\alpha\beta}\circ f_{\alpha}=\rho_{\beta\alpha}\circ f_{\beta}$ for each pair of indices $(\alpha,\beta)$, because $U\cap U_{\alpha}\cap U_{\beta}\neq\emptyset$ since $X$ is irreducible, and we have
$\rho_{U\cap U_{\beta}}^{U_{\beta}}\circ\rho_{\beta\alpha}^{-1}
\circ\rho_{\alpha\beta}\circ(\rho_{U\cap U_{\alpha}}^{U_{\alpha}})^{-1}\circ\rho_{U\cap U_{\alpha}}^{U}=\rho_{U\cap U_{\beta}}^{U_{\beta}}\circ(\rho_{U\cap U_{\alpha}\cap U_{\beta}}^{U_{\beta}})^{-1}\circ\rho_{U\cap U_{\alpha}\cap U_{\beta}}^{U_{\alpha}\cap U_{\beta}}\circ\rho_{\alpha\beta}\circ(\rho_{U\cap U_{\alpha}}^{U_{\alpha}})^{-1}\circ\rho_{U\cap U_{\alpha}}^{U}=(\rho_{U\cap U_{\alpha}\cap U_{\beta}}^{U\cap U_{\beta}})^{-1}\circ\rho_{U\cap U_{\alpha}\cap U_{\beta}}^{U_{\alpha}}\circ(\rho_{U\cap U_{\alpha}}^{U_{\alpha}})^{-1}\circ\rho_{U\cap U_{\alpha}}^{U}=(\rho_{U\cap U_{\alpha}\cap U_{\beta}}^{U\cap U_{\beta}})^{-1}\circ\rho_{U\cap U_{\alpha}\cap U_{\beta}}^{U}=\rho_{U\cap U_{\beta}}^{U}$. Since $\mathscr{F}$ is a sheaf therefore there exists a (unique) morphism $f:\mathscr{F}(U)\rightarrow\mathscr{F}(X)$ so that for each $\alpha$, $f_{\alpha}=\rho_{\alpha}\circ f$. Now it is easy to see that $\rho_{U}^{X}\circ f$ is the identity morphism, because for each $\alpha$ one has $\rho_{U\cap U_{\alpha}}^{U}\circ\rho_{U}^{X}\circ f=\rho_{U\cap U_{\alpha}}^{X}\circ f=\rho_{U\cap U_{\alpha}}^{U_{\alpha}}\circ\rho_{\alpha}\circ f=\rho_{U\cap U_{\alpha}}^{U_{\alpha}}\circ(\rho_{U\cap U_{\alpha}}^{U_{\alpha}})^{-1}\circ\rho_{U\cap U_{\alpha}}^{U}=\rho_{U\cap U_{\alpha}}^{U}$. Moreover, $f\circ\rho_{U}^{X}$ is the identity morphism, because for each $\alpha$ one has $\rho_{\alpha}\circ f\circ\rho_{U}^{X}=f_{\alpha}\circ\rho_{U}^{X}=(\rho_{U\cap U_{\alpha}}^{U_{\alpha}})^{-1}\circ\rho_{U\cap U_{\alpha}}^{U}\circ\rho_{U}^{X}=(\rho_{U\cap U_{\alpha}}^{U_{\alpha}})^{-1}\circ\rho_{U\cap U_{\alpha}}^{X}=\rho_{\alpha}$. \\
\end{bk}

\section{Inverse images of presheaves of groups or rings}

\begin{bk}\label{I: 3.7.1} We show that when $\emph{\textbf{K}}$ is the category of sets then the inverse image of any presheaf $\mathscr{G}$ (over $Y$ with values in $\emph{\textbf{K}}$) by $\psi:X\rightarrow Y$ exists (for the notations consider \ref{I: 3.5.3}). In fact, for each open subspace $U$ of $X$ we define $\mathscr{G}'(U)$ as a subset of $\prod\limits_{x\in U}\mathscr{G}_{\psi(x)}$  to be the set of all of the sequences $( s(x) )_{x\in U}\in\prod\limits_{x\in U}\mathscr{G}_{\psi(x)}$ (for each $x\in U$, $s(x)\in\mathscr{G}_{\psi(x)}$) so that for each $x\in U$ there exists an open neighborhood $V$ of $\psi(x)$ in Y, an open neighborhood $W$ of $x$ contained in $U\cap\psi^{-1}(V)$ and an element $t\in\mathscr{G}(V)$ such that $s(z)=t_{\psi(z)}$ for each $z\in W$ (recall that $t_{\psi(z)}\in\mathscr{G}_{\psi(z)}$ means that the germ of $t$ at the point $\psi(z)$ ).
For each two open subspaces $U'$ and $U$ of $X$ with $U'\subseteq U$, take $(\rho')_{U'}^{U}:\mathscr{G}'(U)\rightarrow\mathscr{G}'(U')$ to be the usual restriction which is well-defined. One can then observe that the assignment $U\rightsquigarrow\mathscr{G}'(U)$ is actually a sheaf over $X$ with values in
$\emph{\textbf{K}}$. We show that the pair $(\mathscr{G}',\rho_{\mathcal{G}})$ where $\rho_{\mathcal{G}}:\mathscr{G}\rightarrow\psi_{\ast}(\mathscr{G}')$ is the canonical morphism in which for each open subspace $V$ of $Y$ is given by $$(\rho_{\mathcal{G}})_{V}:\mathscr{G}(V)
\rightarrow(\psi_{\ast}(\mathscr{G}'))(V)=\mathscr{G}'(\psi^{-1}(V))$$
$$s\rightsquigarrow(s_{\psi(x)})_{x\in\psi^{-1}(V)},$$ is the inverse image of $\mathscr{G}$ by $\psi$. Let $\mathscr{F}$ be a sheaf over $X$ with values in
$\emph{\textbf{K}}$, by definition it is enough to show that the function $$\Hom_{X}(\mathscr{G}',\mathscr{F})\rightarrow
\Hom_{Y}(\mathscr{G},\psi_{\ast}(\mathscr{F}))$$ $$u\rightsquigarrow\psi_{\ast}(u)\circ\rho_{\mathcal{G}}$$ is bijective. For injectiveness suppose that there exist two morphisms
$u, u'\in\Hom_{X}(\mathscr{G}',\mathscr{F})$ so that $\psi_{\ast}(u)\circ\rho_{\mathcal{G}}=\psi_{\ast}(u')\circ\rho_{\mathcal{G}}$. Let $U$ be an open subspace of $X$ and let $\textbf{s}=(s(x))_{x\in U}\in\mathscr{G}'(U)$ be an arbitrary element, therefore for each $p\in U$ there exists an open subspace $V_{p}$ of $Y$ with $\psi(p)\in V_{p}$, an open subspace $W_{p}$ of $X$ with $p\in W_{p}\subseteq\psi^{-1}(V_{p})\cap U$ and an element $t^{p}\in\mathscr{G}(V_{p})$ so that for each $z\in W_{p}$, $s(z)=(t^{p})_{\psi(z)}$. Since $\mathscr{F}$ is a sheaf and $U=\bigcup\limits_{p\in U}W_{p}$ therefore to prove $u(\textbf{s})=u'(\textbf{s})$ it suffices to show that  $u(\textbf{s})|_{W_{p}}=u'(\textbf{s})|_{W_{p}}$ for each $p\in U$.
But one has $u(\textbf{s})|_{W_{p}}=u((s(x))_{x\in U})|_{W_{p}}=u((s(x))_{x\in U}|_{W_{p}})=u((s(x))_{x\in W_{p}})=u(((t^{p})_{\psi(x)})_{x\in W_{p}})=u_{W_{p}}(\rho_{\mathcal{G}}(t^{p})|_{W_{p}})=
u_{\psi^{-1}(V_{p})}(\rho_{\mathcal{G}}(t^{p}))=
(\psi_{\ast}(u))_{V_{p}}(\rho_{\mathcal{G}}(t^{p}))=
(\psi_{\ast}(u)\circ\rho_{\mathcal{G}})_{V_{p}}(t^{p})=
(\psi_{\ast}(u')\circ\rho_{\mathcal{G}})_{V_{p}}(t^{p})=
u'(\textbf{s})|_{W_{p}}$. For surjectiveness, let $v\in\Hom_{Y}(\mathscr{G},\psi_{\ast}(\mathscr{F}))$ be an arbitrary morphism. Let $U$ be an arbitrary open subspace of $X$ we define the function $u_{U}:\mathscr{G}'(U)\rightarrow\mathscr{F}(U)$ as follows:
for each element $\textbf{s}=(s(x))_{x\in U}\in\mathscr{G}'(U)$,
taking into account the previous notations, then for each $p\in U$ set
$t(p)=v(t^{p})|_{W_{p}}\in\mathscr{F}(W_{p})$, we observe that $t(p)|_{W_{p}\cap W_{q}}=t(q)|_{W_{p}\cap W_{q}}$
for each $p,q\in U$, because for each $z\in W_{p}\cap W_{q}$ one has $(t(p)|_{W_{p}\cap W_{q}})_{z}=(v(t^{p})|_{W_{p}\cap W_{q}})_{z}=(v(t^{p}))_{z}=\psi_{z,\mathscr{F}}((v(t^{p}))_{\psi(z)})=
\psi_{z,\mathscr{F}}(v_{\psi(z)}((t^{p})_{\psi(z)}))=
\psi_{z,\mathscr{F}}(v_{\psi(z)}(s(z)))=
\psi_{z,\mathscr{F}}(v_{\psi(z)}((t^{q})_{\psi(z)}))=(t(q)|_{W_{p}\cap W_{q}})_{z}$. Since $\mathscr{F}$ is a sheaf hence there exists a unique element $t\in\mathscr{F}(U)$ such that for each $p\in U$, $t|_{W_{p}}=t(p)$; we define $u_{U}(\textbf{s})=t$. One can then observe that $u:\mathscr{G}'\rightarrow\mathscr{F}$ is actually a morphism of sheaves and $\psi_{\ast}(u)\circ\rho_{\mathcal{G}}=v$. \\
Therefore we have the inverse image functor $\psi^{\ast}:\psh_{\emph{\textbf{K}}}(Y)\rightarrow
\sh_{\emph{\textbf{K}}}(X)$. Moreover, let $u:\mathscr{G}_{1}\rightarrow\mathscr{G}_{2}$ be an arbitrary morphism of presheaves of sets over $Y$; for each open subspace $U$ of $X$ define the function $\theta_{U}:(\psi^{\ast}(\mathscr{G}_{1}))(U)\rightarrow(\psi^{\ast}
(\mathscr{G}_{2}))(U)$ by $(s(x))_{x\in U}\rightsquigarrow(u_{\psi(x)}(s(x)))_{x\in U}$; we observe that $\theta:\psi^{\ast}(\mathscr{G}_{1})\rightarrow\psi^{\ast}
(\mathscr{G}_{2})$ is actually a morphism of sheaves and $\psi_{\ast}(\theta)\circ\rho_{\mathcal{G}_{1}}=\rho_{\mathcal{G}_{2}}\circ u$, therefore $\psi^{\ast}(u)=\theta$. Furthermore for each open subspace $V$ of $Y$, the inverse image of $\mathscr{G}|_{V}$ by the restriction of $\psi$ to $\psi^{-1}(V)$ is identical to the induced sheaf $\psi^{\ast}(\mathscr{G})|_{\psi^{-1}(V)}$, i.e. $$(\psi|_{\psi^{-1}(V)})^{\ast}(\mathscr{G}|_{V})=
\psi^{\ast}(\mathscr{G})|_{\psi^{-1}(V) } .$$
When $\psi$ is the identity map $1_{X}$, we recover the definition of a sheaf of sets associated to a presheaf (of sets over $X$). The above considerations apply without change when $\emph{\textbf{K}}$ is the category of groups or the category of rings (not necessarily commutative).\\
When $X$ is an arbitrary subspace of $Y$, and $j:X\rightarrow Y$ is the canonical injection, for each sheaf $\mathscr{G}$ over $Y$ with values in $\emph{\textbf{K}}$, the inverse image sheaf $j^{\ast}(\mathscr{G})$, if it exists, is said to be the \emph{induced} sheaf over $X$ by $\mathscr{G}$. For the category of sets (or of groups, or of rings) the induced sheaf over $X$ by every sheaf over $Y$ exists.\\
\end{bk}

\begin{bk}\label{I: 3.7.2} Keeping the notations and hypotheses of \ref{I: 3.5.3}, suppose that $\mathscr{G}$ is a presheaf of  sets (or of groups or of rings) over $Y$, the definition of the sections of $\psi^{\ast}(\mathscr{G})$ shows that (also taking into account \ref{I: 3.4.4}) the composition of the canonical homomorphisms of the fibers $\psi_{x}
\circ\rho_{\psi(x)}=\psi_{x,\psi^{\ast}(\mathscr{G})}
\circ(\rho_{\mathcal{G}})_{\psi(x)}:\mathscr{G}_{\psi(x)}\rightarrow
(\psi^{\ast}(\mathscr{G}))_{x}$ is a functorial isomorphism at $\mathscr{G}$; which permits us to identify these two fibers. In fact, one can easily observe that the rule of the morphism $\psi_{x}\circ\rho_{\psi(x)}:\mathscr{G}_{\psi(x)}\rightarrow
(\psi^{\ast}(\mathscr{G}))_{x}$ is of the form $$\langle V, s\rangle\rightsquigarrow\langle\psi^{-1}(V), (s_{\psi(p)})_{p\in\psi^{-1}(V)}\rangle$$ where $V$ is an open neighborhood of $\psi(x)$ in $Y$ and $s\in\mathscr{G}(V)$, it is obviously bijective; this in particular, when $\emph{\textbf{K}}$ is the category of groups or rings, implies that $\Supp(\psi^{\ast}(\mathscr{G}))=\psi^{-1}(\Supp(\mathscr{G}))$. Moreover for each morphism $u:\mathscr{G}\rightarrow\mathscr{H}$ of such presheaves the following diagram is commutative $$\xymatrix{
\mathscr{G}_{\psi(x)}\ar[r]^{\psi_{x}\circ\rho_{\psi(x)}} \ar[d]^{u_{\psi(x)}} & (\psi^{\ast}(\mathscr{G}))_{x} \ar[d]^{(\psi^{\ast}(u))_{x}} \\ \mathscr{H}_{\psi(x)}\ar[r]^{\psi_{x}\circ\rho_{\psi(x)}} & (\psi^{\ast}(\mathscr{H}))_{x} .} $$
We also observe that for each morphism (of presheaves) $u:\mathscr{G}\rightarrow\psi_{\ast}(\mathscr{F})$ where
$\mathscr{F}$ is a sheaf over $X$ with values in $\emph{\textbf{K}}$, the fiber of the corresponding morphism $u^{\sharp}:\psi^{\ast}(\mathscr{G})\rightarrow\mathscr{F}$ at $x$ is of the form: $$u^{\sharp}_{x}=\psi_{x,\mathscr{F}}\circ u_{\psi(x)}\circ(\psi_{x,\psi^{\ast}(\mathscr{G})}
\circ(\rho_{\mathcal{G}})_{\psi(x)})^{-1}:
(\psi^{\ast}(\mathscr{G}))_{x}\rightarrow\mathscr{F}_{x}.$$ An immediate consequence from this result is that the functor $\psi^{\ast}:\psh_{\Ab}(Y)\rightarrow
\sh_{\Ab}(X)$ is exact where $\Ab$ is the category of abelian groups.
\end{bk}

\section{Sheaves of pseudo-discrete spaces}

In this section, $X$ is a topological space whose topology admits a basis $\mathfrak{B}$ formed by quasi-compact open subspaces.\\

\begin{bk}\label{I: 3.8.1}  Let $\mathscr{F}$ be a sheaf of sets over $X$; if for each open subspace $U$ of $X$ we equip $\mathscr{F}(U)$ to the discrete topology, then $U\rightsquigarrow\mathscr{F}(U)$ is a presheaf of topological spaces over $X$. We shall observe that there is a sheaf of topological spaces (over $X$)  $\mathscr{F}'$ associated to the presheaf $\mathscr{F}$ (\ref{I: 3.5.6}) so that for each quasi-compact open subspace $U$ of $X$ (note that $U$ does not necessarily belong to $\mathfrak{B}$), $\mathscr{F}'(U)$ is the discrete space $\mathscr{F}(U)$. First we show that the presheaf $U\rightsquigarrow\mathscr{F}(U)$ of discrete topological spaces over $\mathfrak{B}$ satisfies in the condition $\mathbf{(F_{0})}$. Therefore let $(U_{\alpha})$ be an open covering for $U\in\mathfrak{B}$ formed by elements of $\mathfrak{B}$ contained in $U$, also let $T$ be an arbitrary topological space, the category of topological spaces is denoted by $\Top$, the function $\Hom_{\Top}(T,\mathscr{F}(U))\rightarrow\prod\limits_{\alpha}
\Hom_{\Top}(T,\mathscr{F}(U_{\alpha}))$ given by $f\rightarrow\rho_{\alpha}\circ f$ is a bijection from the set $\Hom_{\Top}(T,\mathscr{F}(U))$ over the set of families $(f_{\alpha})\in\prod\limits_{\alpha}
\Hom_{\Top}(T,\mathscr{F}(U_{\alpha}))$ so that $\rho_{V}^{U_{\alpha}}\circ f_{\alpha}=\rho_{V}^{U_{\beta}}\circ f_{\beta}$ for each pair of indices $(\alpha,\beta)$ and for each $V\in\mathfrak{B}$ with $V\subseteq U_{\alpha}\cap U_{\beta}$. Since $\mathscr{F}$ is a sheaf of sets the foregoing function is obviously injective; also for each such family $(f_{\alpha})$ there exists a (unique) function $f:T\rightarrow\mathscr{F}(U)$ so that $f_{\alpha}=\rho_{\alpha}\circ f$ for each $\alpha$; the function $f$ is continuous; because since $U$ is quasi-compact therefore a finite subset of indices, say $\{\alpha_{1},...,\alpha_{n}\}$, can be extracted in which $U=\bigcup\limits_{j=1}^{n}U_{\alpha_{j}}$; one can then observe that for each $s\in\mathscr{F}(U)$, $$f^{-1}(s)=\bigcap\limits_{j=1}^{n}f_{\alpha_{j}}^{-1}(s|_{U_{\alpha_{j}}}).$$ Since the topologies of $\mathscr{F}(U)$ and $\mathscr{F}(U_{\alpha_{j}})$, for each $1\leq j\leq n$, are the discrete topology therefore the map $f$ is continuous. The category of topological spaces admits the projective limits therefore by \ref{I: 3.2.2}, $\mathscr{F}'$ (with the notations of \ref{I: 3.2.1}) is a sheaf of topological spaces over $X$ (note that for an arbitrary open subspace $U$ of $X$, by the universal property of the projective limit, there exists a (unique) morphism $(\rho_{\mathcal{F}})_{U}:\mathscr{F}(U)\rightarrow\mathscr{F}'(U)$ so that for each $V\in\mathfrak{B}$ with $V\subseteq U$, $\rho_{V}^{U}=\can_{V}\circ(\rho_{\mathcal{F}})_{U}$ where $\can_{V}:\mathscr{F}'(U)\rightarrow\mathscr{F}(V)$ is the canonical morphism. One can observe that $\rho_{\mathcal{F}}:\mathscr{F}\rightarrow\mathscr{F}'$ is a morphism of presheaves. Since $\mathscr{F}$ is a sheaf of sets therefore the canonical morphism $\rho_{\mathcal{F}}:\mathscr{F}\rightarrow\mathscr{F}'$ is an isomorphism as sheaves of sets, therefore for each open subspace $U$ of $X$, we can identify the underlying set of the topological space $\mathscr{F}'(U)$ as $\mathscr{F}(U)$, but the topology of $\mathscr{F}'(U)$ in general is not necessarily discrete: it is the coarsest topology which making the morphism $(\rho')_{V}^{U}:\mathscr{F}'(U)\rightarrow\mathscr{F}'(V)$ continuous where $V\in\mathfrak{B}$ with $V\subseteq U$); one can also observe that the sheaf of topological spaces  $\mathscr{F}'$ together with the canonical morphism $\rho_{\mathcal{F}}:\mathscr{F}\rightarrow\mathscr{F}'$ is the associated sheaf of the presheaf of topological spaces $\mathscr{F}$, i.e. $\mathscr{F}'=1_{X}^{\ast}(\mathscr{F})$ as sheaf of topological spaces. Because let $\mathscr{G}$ be an arbitrary sheaf of topological spaces over $X$, we show that the function $$\theta:\Hom_{\Top}(\mathscr{F}',\mathscr{G})\rightarrow
\Hom_{\Top}(\mathscr{F},\mathscr{G})$$ $$v\rightsquigarrow v\circ\rho_{\mathcal{F}}$$ is bijective. Since $\mathscr{F}$ is a sheaf of sets thus $\mathscr{F}'=1_{X}^{\ast}(\mathscr{F})$ as sheaf of sets (note that every sheaf of topological spaces is also a sheaf of sets), therefore
the function $$\Hom_{\sets}(\mathscr{F}',\mathscr{G})\rightarrow
\Hom_{\sets}(\mathscr{F},\mathscr{G})$$ $$v\rightsquigarrow v\circ\rho_{\mathcal{F}}$$ is bijective. Since $\Hom_{\Top}(\mathscr{F}',\mathscr{G})\subseteq
\Hom_{\sets}(\mathscr{F}',\mathscr{G})$, hence the function $\theta$ is injective, moreover for morphism $u:\mathscr{F}\rightarrow\mathscr{G}$ as presheaves of topological spaces there exists a (unique) morphism $\eta:\mathscr{F}'\rightarrow\mathscr{G}$ as sheaves of sets so that $u=\eta\circ\rho_{\mathcal{F}}$; to conclude the assertion we should show that $\eta:\mathscr{F}'\rightarrow\mathscr{G}$ is a morphism as sheaves of topological spaces. In other words, for each open subspace $U$ of $X$ we should show that $\eta_{U}:\mathscr{F}'(U)\rightarrow\mathscr{G}(U)$ is continuous. Let $U=\bigcup\limits_{\alpha}V_{\alpha}$ be an open covering formed by elements of $\mathfrak{B}$, then for each $\alpha$, take $f_{\alpha}=\eta_{V_{\alpha}}\circ\rho'_{\alpha}:
\mathscr{F}'(U)\rightarrow\mathscr{G}(V_{\alpha})$, note that for each $\alpha$ the morphism $\eta_{V_{\alpha}}:\mathscr{F}'(V_{\alpha})
\rightarrow\mathscr{G}(V_{\alpha})$ is continuous because $\mathscr{F}'(V_{\alpha})$ is canonically isomorphic to $\mathscr{F}(V_{\alpha})$ as topological spaces (\ref{I: 3.2.1}), therefore
the topology of
$\mathscr{F}'(V_{\alpha})$ is discrete. The family of continuous maps $(f_{\alpha})\in\prod\limits_{\alpha}
\Hom_{\Top}(\mathscr{F}'(U),\mathscr{G}(V_{\alpha}))$ satisfies in the condition $\omega_{\alpha\beta}\circ f_{\alpha}=\omega_{\beta\alpha}\circ f_{\beta}$ for each pair of indices $(\alpha,\beta)$ where $\omega_{\alpha\beta}:
\mathscr{G}(V_{\alpha})\rightarrow
\mathscr{G}(V_{\alpha}\cap V_{\beta})$ is the restriction morphism; since $\mathscr{G}$ is a sheaf of topological spaces hence there exists a (unique) continuous map $f:\mathscr{F}'(U)\rightarrow\mathscr{G}(U)$ so that $f_{\alpha}=\omega_{\alpha}\circ f$ for each $\alpha$. Also, $\mathscr{G}$ is a sheaf of sets therefore the function $$\Hom_{\sets}(\mathscr{F}'(U),\mathscr{G}(U))\rightarrow
\prod\limits_{\alpha}
\Hom_{\sets}(\mathscr{F}'(U),\mathscr{G}(V_{\alpha}))$$  $$g\rightsquigarrow(\omega_{\alpha}\circ g)$$ is in particular injective, therefore $\eta_{U}=f$ and so $\eta_{U}$ is continuous (note that since for each open subspace $U$ of $X$ the topology of $\mathscr{F}(U)$ is discrete therefore $\Hom_{\Top}(\mathscr{F},\mathscr{G})=
\Hom_{\sets}(\mathscr{F},\mathscr{G})$, and so by the above correspondence  $\Hom_{\Top}(\mathscr{F}',\mathscr{G})=
\Hom_{\sets}(\mathscr{F}',\mathscr{G})$ ).\\
Now let $U$ be a quasi-compact open subspace of $X$, there exists a finite subset, say $\{U_{1},...,U_{n}\}$, in $\mathfrak{B}$ so that $U=\bigcup\limits_{j=1}^{n}U_{j}$; for each $s\in\mathscr{F}'(U)$, as before, one can observe that $\bigcap\limits_{j=1}^{n}(\rho')_{j}^{-1}(s_{j})=\{s\}$ where $s_{j}=\rho_{j}'(s)\in\mathscr{F}'(U_{j})$; but for each $V\in\mathfrak{B}$, the topology of $\mathscr{F}'(V)$ is discrete, therefore the topology of $\mathscr{F}'(U)$ is so.\\
The preceding considerations can be applied without any change to sheaves of groups or rings (not necessarily commutative), and their associated sheaves respectively are the sheaves of topological groups or topological rings (resp. groups, rings). For shorten, we say that the sheaf $\mathscr{F}'$ is the sheaf of pseudo-discrete spaces (resp. pseudo-discrete groups, pseudo-discrete rings) associated to the sheaf of sets (resp. groups, rings) $\mathscr{F}$.\\
\end{bk}

\begin{bk}\label{I: 3.8.2} Let $\mathscr{F}$ and $\mathscr{G}$ be sheaves of sets (resp. groups, rings) over $X$, $u:\mathscr{F}\rightarrow\mathscr{G}$ a morphism; as we observed in the above (\ref{I: 3.8.1}), we can identify $\mathscr{F}'$ and $\mathscr{G}'$, as sheaves of sets, with $\mathscr{F}$ and $\mathscr{G}$ respectively, and since $\Hom_{\Top}(\mathscr{F}',\mathscr{G}')=
\Hom_{\sets}(\mathscr{F}',\mathscr{G}')$ (\ref{I: 3.8.1});  therefore $u:\mathscr{F}'\rightarrow\mathscr{G}'$ is also a morphism as sheaves of topological spaces (resp. topological groups, topological rings), by abuse of the language $u:\mathscr{F}'\rightarrow\mathscr{G}'$ is also called a continuous homomorphism. More precisely there exists a unique morphism $u':\mathscr{F}'\rightarrow\mathscr{G}'$ as sheaves of topological spaces (resp. topological groups, topological rings) so that $\rho_{\mathcal{G}}\circ u=u'\circ\rho_{\mathcal{F}}$.\\
\end{bk}

\begin{bk}\label{I: 3.8.3} Let $\mathscr{F}$ be sheaf of sets over $X$,  $\mathscr{H}$ a sub-sheaf of $\mathscr{F}$ (i.e. $\mathscr{H}$ is a sheaf of sets so that for each open subspace $U$ of $X$, $\mathscr{H}(U)$ is a subset of $\mathscr{F}(U)$ and for each pair $(U,V)$ of open subspaces of $X$ with $U\subseteq V$ the restriction morphism $\mathscr{H}(V)\rightarrow\mathscr{H}(U)$ is actually the restriction of the morphism $\rho_{U}^{V}:\mathscr{F}(V)\rightarrow\mathscr{F}(U)$ ); let $\mathscr{F}'$ and $\mathscr{H}'$ be sheaves of pseudo-discrete spaces associated to $\mathscr{F}$ and $\mathscr{H}$ respectively; consider the inclusion monomorphism $i:\mathscr{H}\rightarrow\mathscr{F}$ then the corresponding morphism $i':\mathscr{H}'\rightarrow\mathscr{F}'$ (\ref{I: 3.8.2}) is also a monomorphism (equivalently, for each open subspace $U$ of $X$,  $i'_{U}:\mathscr{H}'(U)\rightarrow\mathscr{F}'(U)$ is injective). Then, for each open subspace $U$ of $X$, $\mathscr{H}'(U)$ is a closed subspace of $\mathscr{F}'(U)$ or more precisely, $\Ima(i'_{U})$ is closed:
$$\Ima(i'_{U})=\bigcap\limits_{V}((\rho')_{V}^{U})^{-1}(\Ima(i'_{V}))$$ where $V$ runs through the set of elements of $\mathfrak{B}$ contained in $U$; and $\Ima(i'_{U})$ denotes the image of $\mathscr{H}'(U)$ under the map $i'_{U}$.\\
\end{bk}

\textbf{Acknowledgements.} The author would like to give sincere thanks to the anonymous referee for careful reading of the manuscript.\\

\end{document}